\documentclass[11pt,a4paper]{article}
\usepackage{amsmath, amsfonts, amsthm, amssymb, graphicx, color}
\usepackage{enumerate, hyperref}

\begin{document}

\title
{Regular graphs with few longest cycles}
\author{
{\sc Carol T. ZAMFIRESCU\footnote{Department of Applied Mathematics, Computer Science and Statistics, Ghent University, Krijgslaan 281 - S9, 9000 Ghent, Belgium  and Department of Mathematics, Babe\c{s}-Bolyai University, Cluj-Napoca, Roumania. E-mail address: czamfirescu@gmail.com}}}
\date{}

\maketitle
\begin{center}
\vspace{2mm}
\begin{minipage}{125mm}
{\bf Abstract.} Motivated by work of Haythorpe, Thomassen and the author showed that there exists a positive constant $c$ such that there is an infinite family of 4-regular 4-connected graphs, each containing exactly $c$ hamiltonian cycles. We complement this by proving that the same conclusion holds for planar 4-regular 3-connected graphs, although it does not hold for planar 4-regular 4-connected graphs by a result of Brinkmann and Van Cleemput, and that it holds for 4-regular graphs of connectivity~2 with the constant $144 < c$, which we believe to be minimal among all hamiltonian 4-regular graphs of sufficiently large order. We then disprove a conjecture of Haythorpe by showing that for every non-negative integer $k$ there is a 5-regular graph on $26 + 6k$ vertices with $2^{k+10} \cdot 3^{k+3}$ hamiltonian cycles. 
We prove that for every $d \ge 3$ there is an infinite family of hamiltonian 3-connected graphs with minimum degree $d$, with a bounded number of hamiltonian cycles. It is shown that if a 3-regular graph $G$ has a unique longest cycle $C$, at least two components of $G - E(C)$ have an odd number of vertices on $C$, and that there exist 3-regular graphs with exactly two such components.

\smallskip

{\bf Key words.} Hamiltonian cycle; longest cycle; regular graph; planar graph

\smallskip

\textbf{MSC 2020.} 05C45; 05C07; 05C38; 05C10

\end{minipage}
\end{center}

\vspace{5mm}

\section{Introduction}

A well-known conjecture of Sheehan from 1975 posits that every hamiltonian 4-regular graph has at least two distinct hamiltonian cycles~\cite{Sh75}. As Thomassen points out in~\cite{Th98}, Sheehan's conjecture---combined with results of Smith and Thomason---implies that every hamiltonian regular graph other than a cycle contains at least two hamiltonian cycles. For an overview of results on this conjecture (and its interplay with symmetry) we refer the reader to~\cite{Wa13}.

We will focus here on two relaxations of Sheehan's conjecture. Firstly, for small values of $k$ we describe hamiltonian $k$-regular graphs with few hamiltonian cycles, noting that for $k = 3$ the behaviour is well-understood: already in 1946, Smith showed that every edge in a 3-regular graph is contained in an even number of hamiltonian cycles (see~\cite{Tu46}), so hamiltonian 3-regular graphs contain at least three hamiltonian cycles. By successively expanding vertices into triangles, from any \emph{polyhedral} (i.e.\ planar and 3-connected) graph containing exactly three hamiltonian cycles---such as the tetrahedron---one obtains an infinite family of polyhedral 3-regular graphs containing exactly three hamiltonian cycles. Secondly, we shall discuss regular graphs with a unique longest cycle, with an emphasis on the 3-regular case.

Together with Goedgebeur and Meersman~\cite{GMZ20} we determined, computationally, for every $n \le 21$ the minimum number of hamiltonian cycles in a hamiltonian 4-regular graph on $n$~vertices and confirmed Sheehan's conjecture up to order~21. A deep result of Gir\~{a}o, Kittipassorn, and Narayanan~\cite{GKN19} states that if an $n$-vertex graph with minimum degree at least 3 has a hamiltonian cycle, then it contains another cycle of length at least $n - cn^{4/5}$, where $c > 0$ is an absolute constant. This settles Sheehan's conjecture \emph{asymptotically}.

Entringer and Swart describe an infinite family of graphs containing exactly one hamiltonian cycle, i.e.\ \emph{uniquely hamiltonian} graphs, in which exactly two vertices are 4-valent and all other vertices are cubic~\cite{ES80}; in the same paper, they ask whether Sheehan's conjecture extends to graphs of minimum degree~4. Bondy asked the same question in~\cite[Problem~7.14]{Bo95}. 
In a perhaps surprising denouement, this question turns out to have a negative answer, as proven by Fleischner~\cite{Fl14}. He showed that there exist infinitely many uniquely hamiltonian graphs in which every vertex has degree~4 or 14.

In Section~2, we address two recent conjectures of Haythorpe~\cite{Ha18}. The first of these two conjectures states that the number of hamiltonian cycles in hamiltonian 4-regular graphs increases as a function of the number of vertices. If true, this would verify Sheehan's conjecture. Together with Thomassen, we disproved this conjecture by observing that there are infinitely many 4-regular graphs containing exactly 216 hamiltonian cycles, and by showing that there exists a positive constant $c$ such that there is an infinite family of 4-regular 4-connected graphs containing precisely $c$ hamiltonian cycles~\cite{TZ}. We shall here prove that there is an infinite family of 4-regular graphs of connectivity~2 containing exactly $144 < c$ hamiltonian cycles---we believe 144 to be minimal for large orders. 

The second conjecture posits that for every $k \ge 5$ and every $n \ge k + 3$, all hamiltonian $k$-regular graphs of order~$n$ have at least
$$f(n,k) := (k - 1)^2 [(k - 2)!]^{\frac{n}{k+1}}$$
hamiltonian cycles. We shall prove that this conjecture is invalid by providing, for infinitely many $n$ and every $k \in \{ 5, 6, 7 \}$, hamiltonian $k$-regular graphs on $n$ vertices with fewer than $f(n,k)$ hamiltonian cycles.


Our above contributions disproving Haythorpe's conjectures have connectivity~2. Although the 4-connected case was settled in~\cite{TZ}, the connectivity~3 was left open. Furthermore, despite the fact that Sheehan's conjecture is true when restricted to planar graphs as Bondy and Jackson~\cite{BJ98} proved that every uniquely hamiltonian planar graph contains a vertex of degree 2 or 3, the question remained whether Haythorpe's conjecture on 4-regular graphs holds if we restrict it to planar graphs. We answer both of these questions by showing that there exist infinitely many hamiltonian polyhedral 4-regular graphs of connectivity~3 with a bounded number of hamiltonian cycles. 

In a recent breakthrough, Brinkmann and Van Cleemput~\cite{BV21} proved that planar 4-connected graphs contain an at least linear number of hamiltonian cycles. (It is generally believed that the true lower bound is quadratic---this is realised by \emph{double wheels}, i.e.\ the join of a cycle and the complement of $K_2$.) Thus, it is certain that the asymptotic behaviour changes between planar 4-regular 3-connected and planar 4-regular 4-connected graphs. In this context we also mention recent work of Barish and Suyama~\cite{BS20}, who investigate the complexity of counting hamiltonian cycles in planar 4-regular 4-connected graphs.



In Section 3 we present an infinite family of 4-regular graphs, each containing an odd number of hamiltonian cycles. In Section~4 we address the question whether there is a certain threshold at which the minimum degree forces the presence of a superconstant number of hamiltonian cycles, and show that this is not the case. In Section~5 we prove that there exists an infinite family of hamiltonian bipartite 3-regular cyclically 4-edge-connected graphs with a bounded number of hamiltonian cycles. In Section~6 we treat, as an alternative relaxation of Sheehan's conjecture, regular graphs with a unique longest cycle; in particular, it is shown that if a 3-regular graph $G$ has a unique longest cycle $C$, at least two components of $G - E(C)$ have an odd number of vertices in $C$, and that there exist graphs with exactly two such components. The article concludes with Section~7 in which open problems are discussed.

All graphs in this article are assumed to be connected, unless explicitly stated otherwise. For a graph $G$, we denote by ${\frak H}(G)$ the set of all hamiltonian cycles of $G$, and put $h(G) := |{\frak H}(G)|$. In this paper, in a non-complete $k$-connected graph $G$ a $k$-vertex set $X$ in $G$ is a \emph{$k$-cut} if $G - X$ is disconnected. Let $G$ be a non-complete graph of connectivity $k$, $X$ a $k$-cut in $G$, and $C$ a component of $G - X$. Then $G[V(C) \cup X]$ is called an \emph{$X$-fragment} of $G$. A path with endvertex $v$ is a \emph{$v$-path}, and a $v$-path with endvertex $w \ne v$ is a \emph{$vw$-path}. For a graph $G$ and adjacent vertices $v$ and $w$ in $G$, we will denote both the edge between $v$ and $w$ and the 2-vertex $vw$-path with $vw$---by context it will always be clear what kind of object we are dealing with. For a possibly disconnected graph $G$ we denote by $\omega(G)$ the number of connected components of $G$. For vertices $v,w \in V(G)$, we denote by $G + vw$ the graph obtained by adding the edge $vw$ to $G$ if $vw \notin E(G)$ and otherwise $G + vw := G$. For a set $S$, we say that $A,B \subset S$ \emph{partition} $S$ if $A \cap B = \emptyset$ and $A \cup B = S$.

\section{On two conjectures of Haythorpe}

In this section we treat the following two conjectures of Haythorpe~\cite{Ha18}.

\bigskip

\noindent \textbf{Conjecture 1} (Conjecture~4.2 in~\cite{Ha18}). For $n \ge 8$, all hamiltonian 4-regular graphs of order~$n$ have at least $9 \cdot 2^{\frac{n+2}{6}}$ hamiltonian cycles.

\bigskip

\noindent \textbf{Conjecture 2} (Conjecture~3.1 in~\cite{Ha18}). For $k \ge 5$ and $n \ge k + 3$, all hamiltonian $k$-regular graphs of order~$n$ have at least $f(n,k) := (k - 1)^2 [(k - 2)!]^{\frac{n}{k+1}}$ hamiltonian cycles.

\bigskip

We note that in~\cite{Ha18} neither conjecture asks for the graphs to be hamiltonian, but since various infinite families of non-hamiltonian $k$-regular graphs are known for every $k \ge 3$---among them the famous 4-regular 4-connected family described by Meredith~\cite{Me73}---we have added the hamiltonicity condition in the conjectures' present formulation. Conjecture~1 was recently shown not to be true by Thomassen and the author~\cite{TZ}. In what follows we first complement this result by describing two infinite families of counterexamples to Conjecture~1: one whose members contain precisely $144$ hamiltonian cycles---a number smaller than the one provided in~\cite{TZ}, and which we believe to be minimal for large orders---, the other consisting of polyhedral graphs (while the examples from~\cite{TZ} are either of connectivity~2 and small but non-zero genus, or 4-connected and of large genus). Thereafter, we disprove Conjecture~2 for $k \in \{ 5, 6, 7 \}$.

\subsection{The 4-regular case}

We require a lemma inspired by a simple observation concerning 3-regular graphs, allowing us to infer from the existence of a graph whose hamiltonian cycles satisfy certain conditions the existence of an infinite family of graphs, each with the same number of hamiltonian cycles as the initial graph. This observation was used in~\cite{GMZ20} to show that there exists an infinite family of planar 3-regular cyclically 4-edge-connected graphs, each with exactly four hamiltonian cycles. Moreover, the zig-zag idea we shall use also appears in~\cite{TZ}. 

\newpage

Let $G$ be a hamiltonian graph and $C = abcd$ an induced $4$-cycle in $G$. We call $C$ \emph{good} if

\smallskip

\noindent (i) there is no hamiltonian $ac$-path in $G - b - d$ and no hamiltonian $bd$-path in $G - a - c$;

\smallskip

\noindent (ii) there is no hamiltonian
\begin{itemize}
\itemsep0em
\item $ab$-path in $G - d$,
\item $ad$-path neither in $G - b$ nor in $G - c$,
\item $bc$-path neither in $G - a$ nor in $G - d$,
\item $cd$-path in $G - b$;
\end{itemize}

\smallskip

\noindent (iii) there is no hamiltonian $ac$-path in $G$ and no hamiltonian $bd$-path in $G$;

\smallskip

\noindent (iv) $\{ ad, bc \} \subseteq \bigcap\limits_{{\frak h} \in {\frak H}(G)} E({\frak h})$; and

\smallskip

\noindent (v) the degree of $a$ and $c$ is at least~4, and the degree of $b$ and $d$ is at least~3.

\bigskip

\noindent Throughout the paper, we will refer to these requirements as Conditions~(i)--(v).

\bigskip

\noindent \textbf{Lemma 1.} \emph{Let $G$ be a graph containing a good $4$-cycle $C = abcd$ and such that every vertex in $V(G) \setminus V(C)$ has degree at least~$4$. Then there exists for every positive integer $k$ a graph of minimum degree~$4$, order $|V(G)| + 2k$, and containing exactly $h(G)$ hamiltonian cycles. Denote this infinite family by ${\cal G}$. If every vertex in $V(G) \setminus \{ b, d \}$ is $4$-valent, and $b$ and $d$ are cubic, then every member of ${\cal G}$ is $4$-regular. If $G$ is plane and $C$ a facial cycle, then every member of ${\cal G}$ is planar. For every positive integer $j \le 4$, if $G$ is $j$-connected, then so is every member of ${\cal G}$.}

\bigskip

\noindent \emph{Proof.} We apply the operation illustrated in Fig.~1 to $C$ and thus obtain from $G$ an infinite family of graphs, among which $G'$ is chosen arbitrarily, but is fixed throughout the proof.

\begin{center}
\includegraphics[height=32mm]{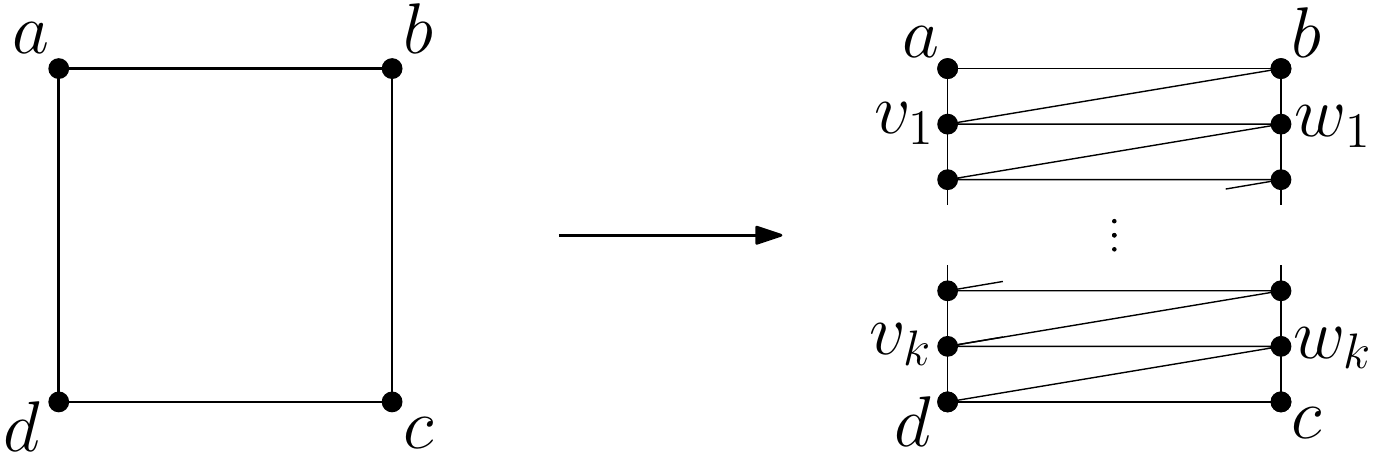}\\[1mm]
Fig.~1: An operation which allows us to construct infinitely many 4-regular graphs with a bounded number of hamiltonian cycles.
\end{center}

We call the vertices $v_1, \ldots, v_k, w_1, \ldots, w_k$ (as defined in Fig.~1) \emph{new} and denote the set of all new vertices by $N$. Since each new vertex is $4$-valent, $G'$ has minimum degree~4, and if every vertex in $V(G) \setminus \{ b, d \}$ is $4$-valent, and $b$ and $d$ are cubic, then $G'$ is $4$-regular. Throughout the proof we see $G - ad - bc$ as a subgraph of $G'$. Put $C' := G'[V(C) \cup N]$. We now show that $h(G) = h(G')$.

Consider distinct ${\frak h}_1, {\frak h}_2 \in {\frak H}(G)$. Let $i \in \{ 1, 2 \}$. By Condition~(iv) we have that $ad, bc \in E({\frak h}_i)$. We now see ${\frak h}_i - ad - bc$ as lying in $G'$. We add to ${\frak h}_i - ad - bc$ the paths $P := av_1 \ldots v_kd$ and $Q := bw_1 \ldots w_kc$, and obtain a hamiltonian cycle ${\frak h}'_i$. We have that ${\frak h}'_1 \ne {\frak h}'_2$ since ${\frak h}'_1 - N = {\frak h}_1 - ad - bc \ne {\frak h}_2 - ad - bc = {\frak h}'_2 - N$. Thus $h(G') \ge h(G)$. Next we show that every hamiltonian cycle in $G'$ can be reduced to a hamiltonian cycle in $G$ and that no two distinct hamiltonian cycles in $G'$ yield the same hamiltonian cycle in $G$. From this it follows that $h(G) \ge h(G')$ and thus $h(G) = h(G')$. Let ${\frak h}'$ be a hamiltonian cycle in $G'$. We consider $S := {\frak h}' \cap C'$, a graph which might be disconnected, and investigate the following three cases, noting that $\omega(S) \le 3$ since $k \ge 1$; we shall treat only essentially different subcases.

\smallskip

\noindent \textsc{Case} 1. $\omega(S) = 1$. Subcase 1.1. $S$ is an $ab$-path spanning $C'$. We consider ${\frak h}' - N + ad + bc$ and obtain a hamiltonian cycle ${\frak h}$ in $G$. As $C$ is chordless, $C'$ contains a unique hamiltonian $ab$-path, so distinct hamiltonian cycles in $G'$ yield distinct hamiltonian cycles in $G$. Subcase~1.2. $S$ is an $ac$-path spanning $C'$. Then ${\frak h}' - N - b - d$ is a hamiltonian $ac$-path in $G - b - d$, contradicting Condition~(i). Subcase~1.3. $S$ is an $ad$-path spanning $C'$. Then ${\frak h}' - N + ab + bc$ is a hamiltonian cycle in $G$ whose intersection with $C$ does not contain $ad$, contradicting Condition~(iv). Subcase~1.4. $S$ is a $bd$-path spanning $C'$. Then ${\frak h}' - N - a - c$ is a hamiltonian $bd$-path in $G - a - c$, contradicting Condition~(i).

\smallskip

\noindent \textsc{Case} 2. $\omega(S) = 2$. Subcase~2.1. $S$ is the disjoint union of an $ab$-path spanning $C' - c$ and $(\{ c \}, \emptyset)$. Then ${\frak h}' - N - d$ is a hamiltonian $ab$-path in $G - d$, contradicting Condition~(ii). Subcase~2.2. $S$ consists of an $ab$-path and a $cd$-path which partition $V(C')$. Then ${\frak h}' - N + ab + cd$ is a hamiltonian cycle in $G$ not containing $ad$, which contradicts Condition~(iv). Subcase~2.3. $S$ is the disjoint union of an $ac$-path spanning $C' - b$ and $(\{ b \}, \emptyset)$. Then ${\frak h}' - N + ad + cd$ is a hamiltonian cycle in $G$ not containing $bc$, which contradicts Condition~(iv). Subcase~2.4. $S$ is the disjoint union of an $ad$-path spanning $C' - b$ and $(\{ b \}, \emptyset)$. Then ${\frak h}' - N - c$ is a hamiltonian $ad$-path in $G - c$, which contradicts Condition~(ii). Subcase~2.5. $S$ is the disjoint union of an $ad$-path spanning $C' - c$ and $(\{ c \}, \emptyset)$. Then ${\frak h}' - N - b$ is a hamiltonian $ad$-path in $G - b$, which contradicts Condition~(ii). Subcase~2.6. $S = P \cup Q$. By considering ${\frak h}' - N + ad + bc$ we obtain a hamiltonian cycle ${\frak h}$ in $G$ such that $E({\frak h}) \cap E(C) = \{ ad, bc \}$. Once more it is clear that, in this situation, distinct hamiltonian cycles in $G'$ yield distinct hamiltonian cycles in $G$. Subcase 2.7. $S$ is the disjoint union of a $bd$-path spanning $C' - a$ and $(\{ a \}, \emptyset)$. Then ${\frak h}' - N + bc + cd$ is a hamiltonian cycle in $G$ not containing $ad$, which contradicts Condition~(iv).

\smallskip

\noindent \textsc{Case} 3. $\omega(S) = 3$. Then $S$ consists of the disjoint union of two isolated vertices and an $xy$-path with $x,y \in \{ a, b, c, d \}$, which together span $C'$. In this situation ${\frak h}'$ either yields a hamiltonian cycle in $G$ whose intersection with $C$ contains only one edge (this occurs if $(x,y) \in \{ (a,b), (a,d), (b,c), (c,d) \}$), in which case we contradict Condition~(iv); or a hamiltonian $ac$-path in $G$ or a hamiltonian $bd$-path in $G$, contradicting Condition~(iii).

\smallskip

The verification of the fact that every member of ${\cal G}$ inherits from $G$ $j$-connectedness for every positive integer $j \le 4$ and planarity if $C$ is a facial cycle is straightforward and left to the reader. \hfill $\Box$

\bigskip

We note that the operation from Fig.~1 may increase the connectivity of the graph it is applied to. We now use Lemma~1 to describe two infinite families of counterexamples to Conjecture~1.

\bigskip

\noindent \textbf{Theorem 1.} \emph{For every integer $n \ge 19$, there exists a $4$-regular graph on $n$~vertices with exactly $144$~hamiltonian cycles. Furthermore, for every integer $n' \ge 22$, there exists a planar $4$-regular graph on $n'$~vertices with exactly $320$~hamiltonian cycles.}

\bigskip

\noindent \emph{Proof.} Consider the graph $G$ shown in Fig.~2 as well as $G' := G - x + dx_1 + cx_2.$ Henceforth, for all objects that are not explicitly defined in the text we refer to Fig.~2. We shall prove that $G$ contains exactly $144$ hamiltonian cycles. The graph $G'$ contains the same number of hamiltonian cycles as $G$; the arguments are very similar to the ones used for $G$ and will therefore be omitted.

Let ${\frak h}$ be an arbitrary but fixed hamiltonian cycle in $G$. The sets $\{ v, w \}$ and $\{ v', w' \}$ are 2-cuts in $G$. Consider, in $G$, the 5-vertex $\{ v, w \}$-fragment $F$ and the 5-vertex $\{ v', w' \}$-fragment $F'$. $F$ ($F'$) contains exactly six hamiltonian $vw$-paths (hamiltonian $v'w'$-paths). Another consequence of the fact that $\{ v, w \}$ and $\{ v', w' \}$ are 2-cuts in $G$ is that $wyv'$, $va$, and $bw'$ are subpaths of ${\frak h}$. Since $wyv' \subset {\frak h}$ we have that $ay, x_4y \notin E({\frak h})$. Therefore ${\frak h}$ must contain $adx$ and $bcx_3$ as subpaths. In $G[\{x, x_1, x_2, x_3, x_4\}]$ there are exactly four hamiltonian $xx_3$-paths. This yields $h(G) = 6^2 \cdot 4$.

\begin{center}
\includegraphics[height=65mm]{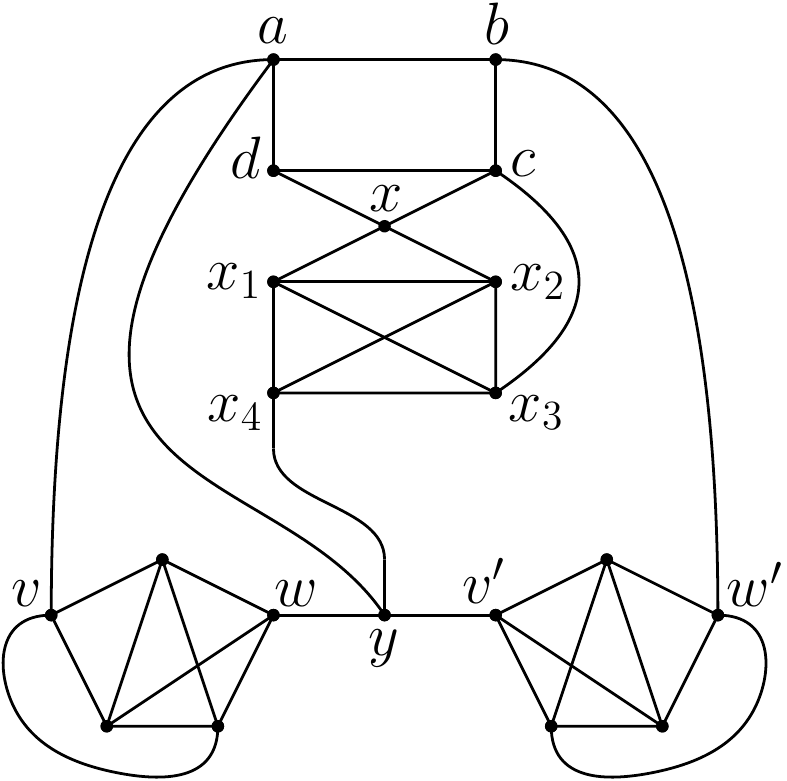}\\
Fig.~2: A graph containing exactly $144$ hamiltonian cycles.\\
As explained in the proof of Theorem~1, the edges $ab$, $ay$, $cd$, $cx$, and $x_4y$ lie on no hamiltonian cycle of the graph, and all other edges lie on some hamiltonian cycle.
\end{center}

We now show that $abcd$ is good in $G$. In $G$, every vertex is $4$-valent with the exception of $b$ and $d$, which are cubic, so Condition~(v) is satisfied. By above arguments we have established that the edges $ad$ and $bc$ lie in every hamiltonian cycle in $G$, so Condition~(iv) is satisfied. In both $G - a - c$ and $G - b - d$ the vertex $y$ is a cut-vertex whose removal yields three components. Hence, there can be no hamiltonian $ac$-path in $G - b - d$ and no hamiltonian $bd$-path in $G - a - c$, so Condition~(i) is satisfied.

In $G - d$, the set $\{ b,y \}$ is a 2-cut and there are three pairwise distinct $\{ b, y \}$-fragments, so there is no hamiltonian $ab$-path in $G - d$ and no hamiltonian $bc$-path in $G - d$. With a similar reasoning one infers that there is no hamiltonian $ad$-path in $G - c$. The graph $G - b$ has connectivity 1 and the cut-vertex $\{ y \}$. Since, in $G - b$, there exists a $\{ y \}$-fragment containing neither $a$ nor $c$ nor $d$, $G - b$ contains no hamiltonian $ad$-path and no hamiltonian $cd$-path. Applying an analogous argument we obtain that $G - a$ does not contain a hamiltonian $bc$-path. Thus, Condition~(ii) is satisfied.

Suppose $G$ contains a hamiltonian $ac$-path ${\frak p}$. Since $\{ a, y \}$ is a 2-cut in $G$, from $a$ the path ${\frak p}$ must visit $v$, traverse $F$, and continue to $y$. As $\{ b, y \}$ is a 2-cut in $G$, from $y$ the path ${\frak p}$ must continue to $v'$, traverse $F'$, and reach $b$. From $b$, however, the edge $ba$ cannot be used and neither can the edge $bc$ be used, as this would contradict the fact that ${\frak p}$ spans $G$. Thus, a contradiction is obtained implying that there is no hamiltonian $ac$-path in $G$. The proof that there is no hamiltonian $bd$-path in $G$ is very similar and left to the reader. Hence, Condition~(iii) is satisfied.

We have proven that $abcd$ is good in $G$. We apply Lemma~1 and obtain for every even integer $n \ge 20$ a 4-regular graph on $n$ vertices containing exactly $h(G) = 144$ hamiltonian cycles. By applying the same reasoning to the 19-vertex graph $G'$ we have completed the proof of the first statement.

For the theorem's second statement, replace each crossing in Fig.~2 by a vertex. \hfill $\Box$ 

\bigskip

Below, we recall information from Table~3 in~\cite{GMZ20}. It implies that Conjecture~1 holds for all $n \le 21$.

\newpage

\begin{center}
\begin{tabular}{ l c c c c c c c c c c c }
Order & 5 & 6 & 7 & 8 & 9,\,10 & 11 & 12 & 13,\,14,\,15,\,16 & 17 & 18 & 19,\,20,\,21 \\
\hline
Nr.~ham.~cycles & 12 & 16 & 23 & 29 & 36 & 48 & 60 & 72 & 96 & 108 & 144
\end{tabular}\\[1mm]
Table 1: For every $n \in \{ 5, \ldots, 21 \}$ the minimum number of hamiltonian cycles occurring in a hamiltonian 4-regular graph on $n$ vertices is given.\\[1mm]

\end{center}

On the other hand, by Theorem~1 and as $9 \cdot 2^{\frac{n+2}{6}} > 144$ for all $n \ge 23$ and $9 \cdot 2^{\frac{n+2}{6}} > 320$ for all $n \ge 29$, we have:

\bigskip

\noindent \textbf{Corollary 1.} \emph{Conjecture 1 does not hold for any $n \ge 23$. Restricted to planar graphs, Conjecture 1 does not hold for any $n \ge 29$.}

\bigskip

We mention in this context a problem of Thomassen raised in~\cite{Th96} and asking whether there exists a 4-regular \emph{bipartite} hamiltonian graph with more than $10^{10}$ vertices and less than 100 hamiltonian cycles.

For a graph $G$ and $M \subset E(G)$ denote by $h(G,M)$ the number of hamiltonian cycles ${\frak h}$ in $G$ satisfying $M \subset E({\frak h})$. Abusing notation, whenever $M$ consists of a single edge $e$ we write $h(G, e)$. We now briefly comment upon which graphs can occur as induced subgraphs in infinite families of hamiltonian 4-regular graphs with a bounded number of hamiltonian cycles. Let $G$ be a $4$-regular graph and $vw$ an edge in $G$. By replacing in the construction from Fig.~2 the 5-vertex $\{ v, w \}$-fragment with $G - vw$, we obtain the following observation.

\bigskip

\noindent \textbf{Proposition 1.} \emph{For any $4$-regular graph $G$ and any edge $e$ in $G$ there exists an infinite family of $4$-regular graphs, each containing $G - e$ as an induced subgraph and exactly $24 \cdot h(G,e)$ hamiltonian cycles.}

\bigskip



2-cuts play a crucial role in Theorem~1. In~\cite{TZ} it is proven that 4-connected counterexamples to Conjecture~1 exist. We now discuss the connectivity 3 case. We shall make use of a fragment which, in its planar version, already appears in the recent description of the smallest known non-hamiltonian polyhedral $4$-regular graph~\cite{VZ18}.

\bigskip

\noindent \textbf{Theorem 2.} \emph{There are infinitely many $4$-regular graphs of connectivity~$3$, each containing exactly $2^{12} \cdot 5 = 20\,480$ hamiltonian cycles. Furthermore, there are infinitely many polyhedral $4$-regular graphs, each containing exactly $2 \cdot 5^5 \cdot 13 = 81\,250$ hamiltonian cycles.}

\bigskip

\noindent \emph{Proof.} Consider the graph $G$ shown in Fig.~3. Henceforth, for all objects that are not explicitly defined in the text we refer to Fig.~3. The graph $G$ has connectivity~3 and all of its vertices are $4$-valent with the exception of $b$ and $d$, which are cubic. We first show that $G$ contains exactly $2^{12} \cdot 5$ hamiltonian cycles.

Consider the 3-cut $X := \{ x, y, z \}$ of $G$. There are two $X$-fragments which we denote by $F$ and $F'$, where $F'$ shall contain $a$. It is straightforward but crucial to note that for any distinct $i,j \in X$, every $ij$-path in $F$ that is a subpath of a hamiltonian cycle of $G$ must contain $k \in X \setminus \{ i, j \}$. There are three essentially different such paths which we will call, noting their endpoints as indices, ${\frak p}_{xy}, {\frak p}_{xz}, {\frak p}_{yz}$. We remark that for a fixed pair of endpoints, we have $2^7$ ways to traverse $F$ from one endpoint to another whilst visiting all of $F$'s vertices. From the above observation we can also infer the following Claim, wherein $L$ and $R$ are induced subgraphs of $G$ isomorphic to $K_4$ located as shown in Fig.~3.

\smallskip

\noindent \textsc{Claim.} \emph{Let $F_L$ be the $\{ a, y, z \}$-fragment of $G$ containing $L$. Consider a subgraph of $G$ containing an $ay$-path $P$ such that $P$ contains an $az$-path spanning $F_L - y$. Then $P$ cannot contain all vertices of $F$ and $R$.}

\smallskip

Let ${\frak h}$ be a hamiltonian cycle in $G$.

\smallskip

\noindent \textsc{Case 1.} ${\frak p}_{xy} \subset {\frak h}$. Since $\{ a, y, z \}$ is a 3-cut in $G$, from $y$ the cycle ${\frak h}$ must enter $L$---in particular, $yx_4 \notin E({\frak h})$---and exit it towards $a$, and there are four ways to do so. Analogously, as $\{ b, x, z \}$ is a 3-cut in $G$, from $x$ the cycle ${\frak h}$ must enter $R$ and exit it towards $b$, and there are four ways to do so. Since ${\frak h}$ is hamiltonian, $ad$ and $bc$ lie in ${\frak h}$ (and $ab$ and $cd$ do not). Up to this point, we have shown that in $G - \{ x_1, x_2, x_3, x_4 \}$ there are exactly $2^{11}$ hamiltonian $cd$-paths containing ${\frak p}_{xy}$.

\begin{center}
\includegraphics[height=82mm]{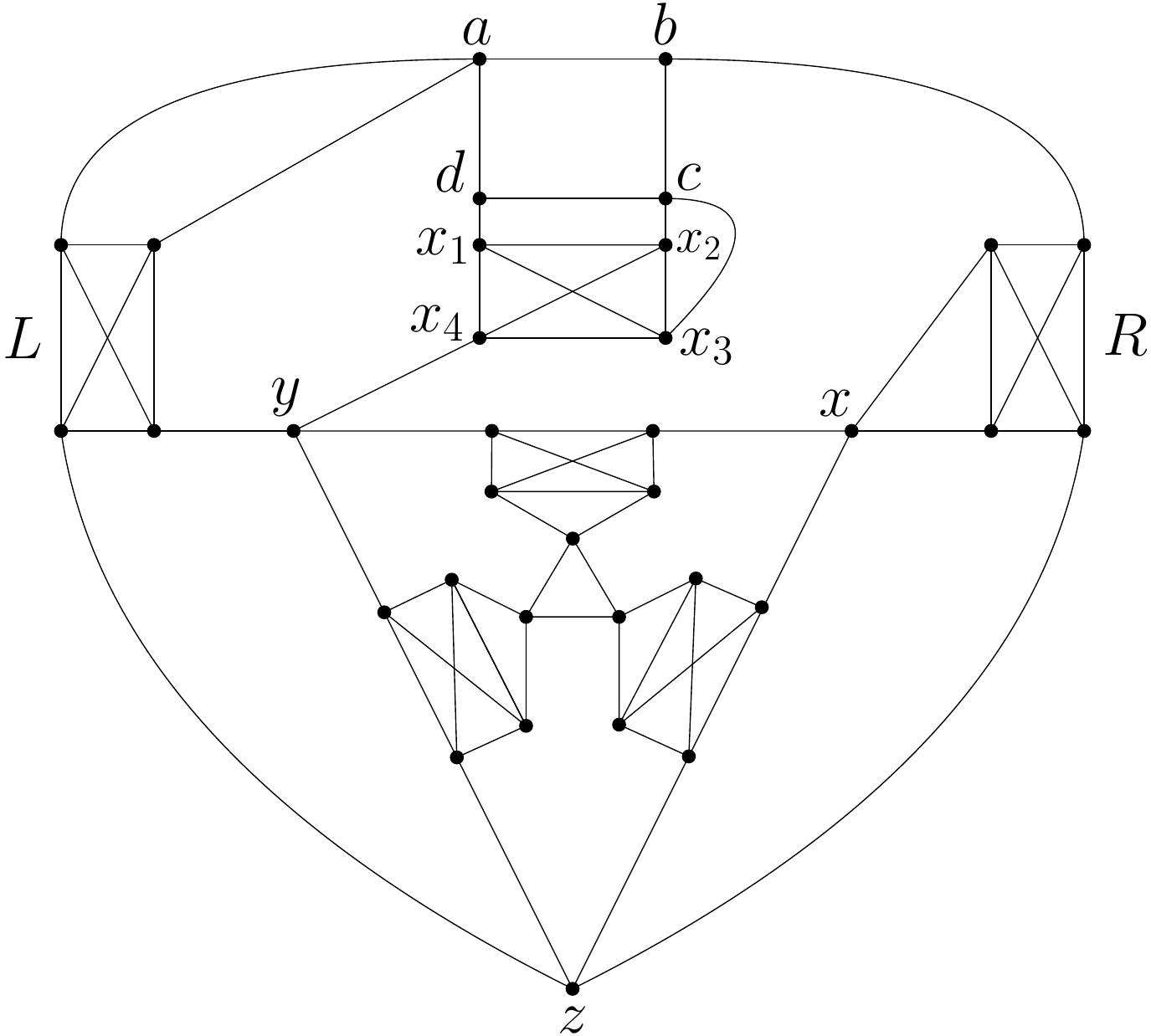}\\
Fig.~3: A graph of connectivity~3 and containing exactly $2^{12} \cdot 5$ hamiltonian cycles.\\
As explained in the proof of Theorem 2, the edges $ab, cd$, and $x_4y$ lie on no hamiltonian
cycle of the graph.
\end{center}

\smallskip

\noindent \textsc{Case 2.} ${\frak p}_{xz} \subset {\frak h}$. As $\{ a, y, z \}$ is a 3-cut in $G$, from $z$ the cycle ${\frak h}$ must enter $L$ and exit it towards $a$ and there are four ways to do so, and since $\{ b, x, z \}$ is a 3-cut in $G$, from $x$ the cycle ${\frak h}$ must enter $R$ and exit it towards $b$; again there are four ways to do so. We have shown that in $G - \{ x_1, x_2, x_3, x_4 \}$ there are exactly $2^{11}$ hamiltonian $cd$-paths containing ${\frak p}_{xz}$.

\smallskip

\noindent \textsc{Case 3.} ${\frak p}_{yz} \subset {\frak h}$. As $\{ b, x, z \}$ is a 3-cut in $G$, from $z$ the cycle ${\frak h}$ must enter $R$ and exit it towards $b$ and there are two ways to do so, and since $\{ a, y, z \}$ is a 3-cut in $G$, from $y$ the cycle ${\frak h}$ must enter $L$ and exit it towards $a$; there are four ways to do so. We have shown that in $G - \{ x_1, x_2, x_3, x_4 \}$ there are exactly $2^{10}$ hamiltonian $cd$-paths containing ${\frak p}_{yz}$.

\smallskip

There are exactly four hamiltonian $cd$-paths in $H := G[\{ c, d, x_1, x_2, x_3, x_4 \}]$. In total we thus obtain $h(G) = 2^{12} \cdot 5$. In order to apply Lemma~1, it remains to show that $abcd$ is good. In $G$, every vertex is $4$-valent with the exception of $b$ and $d$, which are cubic, so Condition~(v) is satisfied. From the above discussion we can infer that the edges $ad$ and $bc$ lie in every hamiltonian cycle of $G$, so Condition~(iv) is also satisfied.

Suppose there exists a hamiltonian $ac$-path ${\frak p}$ in $G - b - d$. The path ${\frak p}$ must exit $H - d$ via $x_4$ towards $y$. We obtain a contradiction to the Claim. The proof that there is no hamiltonian $bd$-path in $G - a - c$ is similar and left to the reader. We have shown that Condition~(i) is satisfied.

Assume there is a hamiltonian $ab$-path ${\frak p}'$ in $G - d$. From $b$, the path ${\frak p}'$ must visit $H - d$ and exit it via $x_4$ towards $y$. Since $ab \notin E({\frak p}')$, by the Claim we obtain a contradiction. The fact that there is no hamiltonian $ad$-path in $G - b$ follows in a very similar way. Since a hamiltonian $ad$-path in $G - c$ cannot use the edge $ad$, such a path cannot exist by the Claim. Suppose there is a hamiltonian $bc$-path ${\frak p}''$ in $G - a$. As $bc \notin E({\frak p}'')$, from $c$ the path ${\frak p}''$ visits $H$ and exits it via $x_4$ towards $y$. From $b$, the path ${\frak p}''$ must visit $R$ and exit it towards $x$ or $z$. In either case, it is now impossible for ${\frak p}''$ to visit all vertices of $L$ and $F$ due to the observation made in this proof's second paragraph, and thus we obtain a contradiction. Similarly, a hamiltonian $bc$-path in $G - d$ cannot contain the edge $bc$, so from $c$ such a path would visit $H - d$ exiting it towards $y$, from where we must visit all vertices of $F$ and reach either $x$ or $z$. But then every vertex in $L$ or $R$ is missed, which yields a contradiction. Finally, a hamiltonian $cd$-path in $G - b$ must contain subpaths between $c$ and $y$ as well as $d$ and $z$. However, in this situation a vertex in $F$ or a vertex in $R$ is not contained in this supposedly hamiltonian path---a contradiction. We have proven that Condition~(ii) is satisfied.

Suppose there is a hamiltonian $ac$-path ${\frak p}_{ac}$ in $G$. If $cd \in E({\frak p}_{ac})$ then ${\frak p}_{ac} \cap H$ is a hamiltonian $cx_4$-path and $x_4y \in E({\frak p}_{ac})$. But then we obtain a contradiction to the Claim. So $cd \notin E({\frak p}_{ac})$, whence, $adx_1 \subset {\frak p}_{ac}$ (since $d$ is cubic). As $c$ is an endpoint of ${\frak p}_{ac}$ we have that ${\frak p}_{ac} \cap (H - c)$ is a hamiltonian $dx_4$-path and $x_4y \in E({\frak p}_{ac})$. Since $a$ is an endpoint of ${\frak p}_{ac}$ and $\{ a, y, z \}$ is a 3-cut of $G$, the path ${\frak p}_{ac}$ must continue from $y$ to $L$, visit all of its vertices, and exit it towards $z$. However, it is now impossible for ${\frak p}_{ac}$ to span $F$, a contradiction. Assume there is a hamiltonian $bd$-path ${\frak p}_{bd}$ in $G$. If $ab \in E({\frak p}_{bd})$, then ${\frak p}_{bd} \cap H$ must be a hamiltonian $dx_4$-path and $x_4y \in E({\frak p}_{bd})$. We obtain a contradiction to the Claim, so $ab \notin E({\frak p}_{bd})$. If ${\frak p}_{bd}$ contains the edge connecting $b$ to $R$, as $d$ is an endpoint of ${\frak p}_{bd}$ it follows that ${\frak p}_{bd}$ contains, as a subpath, a $yz$-path spanning $F_L$, and also the edge $x_4y$. But then ${\frak p}_{bd}$ cannot contain all vertices of $F$, a contradiction. Therefore $bc \in E({\frak p}_{bd})$. Since $d$ is an endpoint of ${\frak p}_{bd}$, ${\frak p}_{bd} \cap (H - d)$ consists of a hamiltonian $cx_4$-path. From $x_4$ the path ${\frak p}_{bd}$ continues to $y$. Furthermore, $ad \in E({\frak p}_{bd})$. Again we obtain a contradiction to the Claim. We have shown that Condition~(iii) is satisfied, as well.

We have proven that $abcd$ is good in $G$. We apply Lemma~1 and thus complete the proof of the first statement. For the theorem's second statement, replace each crossing in Fig.~3 by a vertex. \hfill $\Box$ 



\bigskip



We conclude this subsection with the following consequence of Theorem~2.

\bigskip

\noindent \textbf{Corollary 2.} \emph{Among $3$-connected graphs, Conjecture 1 does not hold for any $n \ge 65$. Restricted to planar $3$-connected graphs, Conjecture~1 does not hold for any $n \ge 77$.}

\bigskip







\subsection{The 5-, 6-, and 7-regular case}

Within this subsection and the Appendix, whenever a figure depicts a graph $G$ with some edges drawn thicker than others, then the set of thick edges shows $M := \bigcup_{{\frak h} \in {\frak H}(G)} E({\frak h})$ and the set of thin edges illustrates $E(G) \setminus M$.

Let $G$ be a graph containing disjoint edges $ac$ and $uw$ which we call \emph{special}. Define the graph $H$ as $G$ in which we subdivide $ac$ with a vertex $b$ and $uw$ with a vertex $v$, and add the edge $bv$. Put 
$h_{00}^{11}(G; ac, uw) := h(H - bc - vw)$, $h_{01}^{01}(G; ac, uw) := h(G - ac, uw)$, $h_{01}^{10}(G; ac, uw) := h(H - bc - uv)$, $h_{10}^{01}(G; ac, uw) := h(H - ab - vw)$, $h_{10}^{10}(G; ac, uw) := h(G - uw, ac)$, $h_{11}^{00}(G; ac, uw) := h(H - ab - uv)$, and $h_{11}^{11}(G; ac, uw) := h(G, \{ ac, uw \})$. We shall abbreviate $h_{ij}^{k\ell}(G; ac, uw)$ to $h_{ij}^{k\ell}$ whenever the choice of $G$, $ac$, and $uw$ are clear from the context. 
The following lemma is inspired by a technique used by Haythorpe~\cite{Ha18}.

\bigskip

\noindent \textbf{Lemma 2.} \emph{Consider a $k$-regular graph $G$ with $k \ge 3$, and a pair of special edges in $G$ such that $G$ does not contain a $2$-factor consisting of exactly two components, each containing a special edge, and $h_{01}^{01} = h_{01}^{10} = h_{10}^{01} = h_{10}^{10} = 0$. Then there exists for every positive integer $\ell$ a $k$-regular graph on $|V(G)| + \ell(k+1)$ vertices containing exactly
$$((k-2)!)^\ell \left( h_{11}^{11} + (k-2)\left(h_{00}^{11} + h_{11}^{00}\right) \right)$$
hamiltonian cycles.}

\bigskip

\noindent \emph{Proof.} Denote by $ac$ and $uw$ the special edges of $G$. On $ac$ we insert vertices $b_1, \ldots, b_\ell$ (in this order, from $a$ towards $c$) and on $uw$ we insert vertices $v_{11}, v_{12}, v_{21}, v_{22}, \ldots, v_{\ell1}v_{\ell2}$ (in this order, from $u$ towards $w$). We have obtained a graph to which we add the pairwise disjoint complete graphs $K^1, \ldots, K^\ell$, each isomorphic to $K_k$. We identify an edge $e_i$ of $K^i$ with $v_{i1}v_{i2}$ and connect $b_i$ to every vertex in $K^i$ that is not a vertex of $e_i$, for all $i \in \{ 1, \ldots, \ell \}$. We call this graph $G'$. The case $k = 5$ and $\ell = 2$ is shown in Fig.~4.

\begin{center}
\includegraphics[height=50mm]{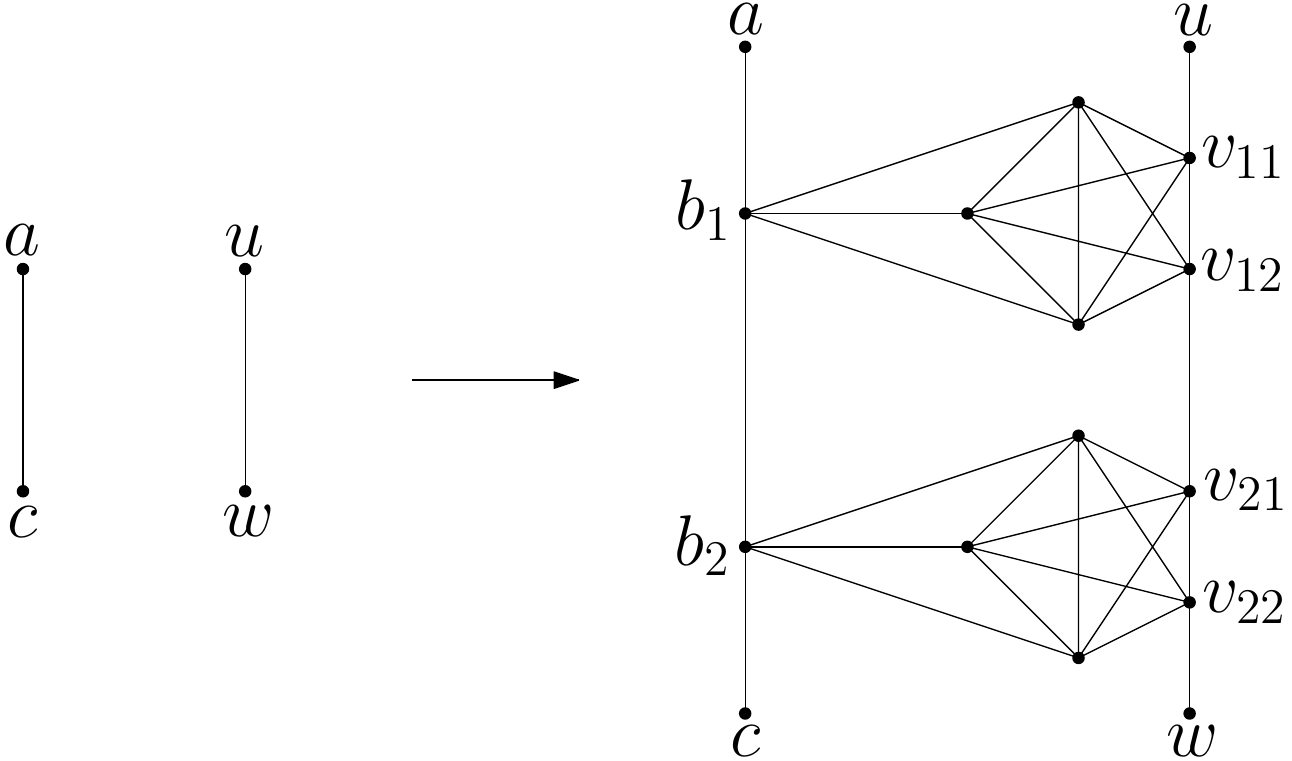}\\
Fig.~4: An illustration of the operation used in the proof of Lemma~2\\
for the case $k = 5$ and $\ell = 2$.
\end{center}

For every $i \in \{ 1, \ldots, \ell \}$, in $K^i$ there are exactly $(k - 2)!$ hamiltonian $v_{i1}v_{i2}$-paths and for every $j \in \{ 1, 2 \}$ there are exactly $(k-2)!(k-2)$ hamiltonian $b_iv_{ij}$-paths in $G'[V(K^i) \cup \{ b_i \}]$. We now treat the five essentially different ways in which a hamiltonian cycle in $G'$ may intersect $J := G'[\{ a, c, u, w, b_1, \ldots, b_\ell \} \cup V(K^1) \cup \ldots \cup V(K^\ell)]$. We call this intersection $S$.

\smallskip

\noindent \textsc{Case 1.} $S$ is either a hamiltonian $uw$-path in $J - a - c$ or a hamiltonian $aw$-path in $J - u - c$ or a hamiltonian $cu$-path in $J - a - w$ or a hamiltonian $ac$-path in $J - u - w$. In these cases we would obtain a contradiction to $h_{01}^{01} = 0$, $h_{01}^{10} = 0$, $h_{10}^{01} = 0$, and $h_{10}^{10} = 0$, respectively.

\smallskip

\noindent \textsc{Case 2.} $S$ consists of an $ac$-path $P_{ac}$ and a $uw$-path $P_{uw}$ whose vertex sets partition $V(J)$. By the counts above, there are exactly $((k-2)!)^\ell h_{11}^{11}$ hamiltonian cycles in $G'$ containing both $P_{ac}$ and $P_{uw}$.

\smallskip

\noindent \textsc{Case 3.} $S$ is a hamiltonian $au$-path ${\frak p}_{au}$ in $J - c - w$. By the counts above there are exactly $(k-2)((k-2)!)^\ell h_{00}^{11}$ hamiltonian cycles in $G'$ containing ${\frak p}_{au}$.

\smallskip

\noindent \textsc{Case 4.} $S$ is a hamiltonian $cw$-path ${\frak p}_{cw}$ in $J - a - u$. By the counts above there are exactly $(k-2)((k-2)!)^\ell h_{11}^{00}$ hamiltonian cycles in $G'$ containing ${\frak p}_{cw}$.

\smallskip

\noindent \textsc{Case 5.} $S$ consists of an $au$-path and a $cw$-path whose vertex sets partition $V(J)$. This immediately implies that in $G$ there is a $2$-factor consisting of exactly two components, each containing a special edge, contradicting one of the hypotheses. \hfill $\Box$

\bigskip

Note that $h(G - ac - uw)$ plays no role in above arguments. We also point out that we could have allowed some or all of $h_{01}^{01}, h_{01}^{10}, h_{10}^{01}, h_{10}^{10}$ to be non-zero in the statement of Lemma~2, as well as permitting a $2$-factor consisting of exactly two components, each containing a special edge, but this would have led to a significantly more technical statement (and proof) without benefit to the application of this lemma, which we now present.

\bigskip

\noindent \textbf{Theorem 3.} \emph{For every non-negative integer $k$ there exists (i) a $5$-regular graph on $26 + 6k$ vertices with exactly $2^{k+10} \cdot 3^{k+3}$ hamiltonian cycles, (ii) a $6$-regular graph on $39 + 7k$ vertices with exactly $2^{3(k+4)} \cdot 3^{k+4} \cdot 5^5$ hamiltonian cycles, and (iii) a $7$-regular graph on $54 + 8k$ vertices with exactly $2^{3(k+8)} \cdot 3^{k+10} \cdot 5^{k+5}$ hamiltonian cycles.}

\bigskip

\noindent \emph{Proof.} (i) Throughout this first part of the proof we refer to Fig.~5 and the therein depicted graph $G$. Let ${\frak h}$ be a hamiltonian cycle in $G$. Each of the bottom three copies of $K_6$ minus an edge can be traversed by ${\frak h}$ in exactly $4!$ ways. By construction, the two vertices marked by white squares in Fig.~5 must be traversed horizontally by ${\frak h}$; we abbreviate this observation by $(\dagger)$. The subgraph of $G$ induced by $a$, $u$, and the four vertices marked with white circular disks in Fig.~5 admits precisely two hamiltonian $au$-paths. So $h(G) = 2 \cdot (4!)^3$.

\begin{center}
\includegraphics[height=65mm]{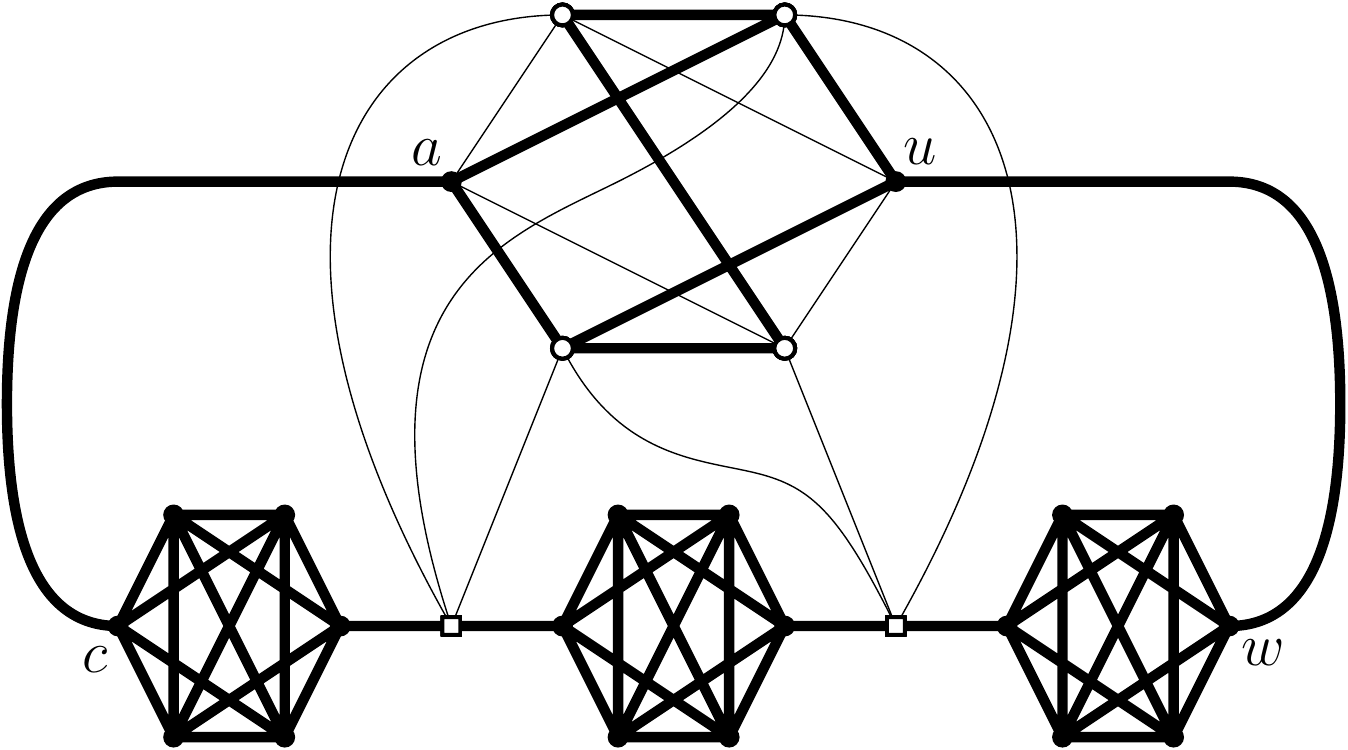}\\
Fig.~5: A 5-regular graph with exactly $2 \cdot (4!)^3$ hamiltonian cycles.
\end{center}

We consider the edges $ac$ and $uw$ to be special. From ($\dagger$) we can infer that $h_{01}^{01} = h_{01}^{10} = h_{10}^{01} = h_{10}^{10} = h_{00}^{11} = h_{11}^{00} = 0$. By construction, it is clear that $G$ does not contain a $2$-factor consisting of exactly two components, each containing a special edge. Thus, by applying Lemma~2, we obtain that there exists for every positive integer $k$ a $5$-regular graph on $26 + 6k$ vertices containing exactly $6^k \cdot h_{11}^{11} = (3!)^k \cdot 2 \cdot (4!)^3$ hamiltonian cycles.

\smallskip

(ii) The proof is very similar to the one given for (i), so we will be succinct. We now refer to Fig.~A from the Appendix and the therein depicted graph $G'$. The subgraph of $G'$ induced by $a$, $u$, and the six vertices marked with white circular disks in Fig.~A admits precisely five hamiltonian $au$-paths. Therefore $h(G') = 5 \cdot (5!)^4$. As above, $G'$ does not contain a $2$-factor consisting of exactly two components, each containing a special edge, so we may apply Lemma~2 and obtain that there exists for every positive integer $k$ a $6$-regular graph on $39 + 7k$ vertices containing exactly $(4!)^k \cdot 5 \cdot (5!)^4$ hamiltonian cycles.

\smallskip

(iii) We apply the same strategy as above but make use of the graph from Fig.~B in the Appendix. \hfill $\Box$

\bigskip











We note that certain variations of the graphs given in Fig.~5 and Figs.~A and B in the Appendix are possible, while maintaining their hamiltonian properties: vertices marked by white squares can be connected to vertices marked by white circular disks or white squares in any way as long as the resulting graph is regular (and, of course, simple).

From Theorem~3, we obtain the following result.

\newpage

\noindent \textbf{Corollary 3.} \emph{Conjecture 2 does not hold for infinitely many $5$-, $6$-, and $7$-regular graphs.}

\bigskip

\noindent \emph{Proof.} The statement follows from Theorem~3~(i), (ii), and (iii), respectively. For the 5-regular case, we note that
$$2^7 \cdot (3!)^{\frac{n-8}{6}} < 2^4 \cdot (3!)^{\frac{n}{6}} = f(n,5)$$
for all $n \ge 26$; for the 6-regular case, we have
$$5^5 \cdot (4!)^{\frac{n-11}{7}} < 5^2 \cdot (4!)^{\frac{n}{7}} = f(n,6)$$
for all $n \ge 39$; and for the 7-regular case, we have
$$2^9 \cdot 3^5 \cdot (5!)^{\frac{n-14}{8}} < 6^2 \cdot (5!)^{\frac{n}{8}} = f(n,7)$$
for all $n \ge 54$. \hfill $\Box$



\section{4-regular graphs with an odd number of hamiltonian cycles}

It is well known that there exist 3-regular graphs containing an even number of hamiltonian cycles and 3-regular graphs containing an odd number of hamiltonian cycles, but that by a theorem of Smith (see~\cite{Tu46}) in a 3-regular graph every edge is traversed by an even number of hamiltonian cycles. For 4-regular graphs, the former holds as well, but the latter does not. Indeed, concerning the former, previous sections discuss various 4-regular graphs containing an even number of hamiltonian cycles, and as can be seen in Table~1, there exist 4-regular graphs on 7 (8; 11; 14; 26) vertices with exactly 23 (29; 145; 323; 25299) hamiltonian cycles. (In fact, the counts up to order~15 were already obtained by Royle and mentioned in~\cite{Ro12}, but not published.) 
We observe that from these examples it is straightforward to infer the existence of infinite families of 4-regular graphs containing an odd number of hamiltonian cycles as follows. This is done by simply forming a ``chain'' of suitable graphs:

\bigskip

\noindent \textbf{Proposition 2.} \emph{For an integer $r \ge 2$, let $S$ be the set of all integers $s$ for which there exists an $r$-regular graph containing an edge which is traversed by exactly $s$ hamiltonian cycles. Then for every set $\{ s_0, \ldots, s_{t-1} \} \subset S$ and any positive integers $k_0, \ldots, k_{t-1}$, there exists an $r$-regular graph containing exactly $\Pi_{i=0}^{t-1} s_i^{k_i}$ hamiltonian cycles.}

\bigskip

\noindent \emph{Proof.} Let $G_i$ be an $r$-regular graph containing an edge $v_iw_i$ traversed by exactly $s_i$ hamiltonian cycles. For every $i \in \{ 0, \ldots, t-1 \}$, consider $k_i$ pairwise disjoint copies of $G_i - v_iw_i$ which we call $H_i^0, \ldots, H_i^{k_i-1}$. We denote the vertices in $H_i^j$ corresponding to $v_i$ ($w_i$) by $v_i^j$ ($w_i^j$). Then the graph
$$\bigcup_{i=0}^{t-1} \left( \bigcup_{j=0}^{k_i-1} H_i^j + w_i^{k_i-1}v_{i+1}^0 + \sum_{\ell = 0}^{k_i - 2}w_i^\ell v_i^{\ell+1} \right)\hspace{-1.2mm},$$
where indices are to be taken mod.~$t$, yields the statement. \hfill $\Box$

\bigskip

\noindent \textbf{Corollary 4.} \emph{There exist infinitely many $4$-regular graphs, each containing an odd number of hamiltonian cycles.}

\bigskip

\noindent \emph{Proof.} The graph from Fig.~6 shows that $11 \in S$, where $S$ is as defined (for $r = 4$) in the statement of Proposition~2. Thus, by applying Proposition~2 we can infer the result. \hfill $\Box$

\begin{center}
\includegraphics[height=45mm]{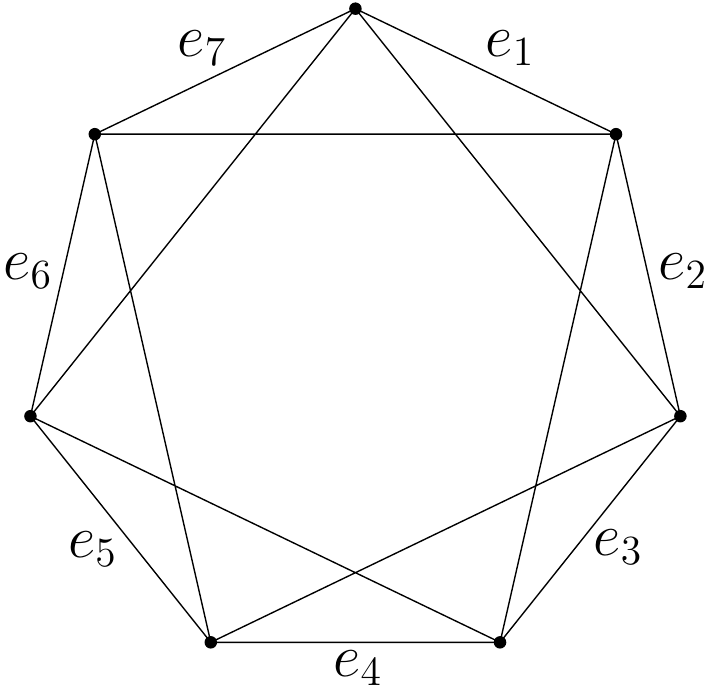}\\
Fig.~6: The antihole on seven vertices is a 4-regular graph containing exactly 23 hamiltonian cycles. Each edge $e_i$ is traversed by exactly 11 hamiltonian cycles, and every other edge is traversed by exactly 12 hamiltonian cycles.
\end{center}

\section{No degree threshold for a superconstant number of hamiltonian cycles}

As mentioned in the introduction, Entringer and Swart, and independently Bondy, asked whether Sheehan's conjecture is true for graphs of minimum degree~4. It is not, as 
Fleischner proved that there exist infinitely many graphs containing exactly one hamiltonian cycle in which every vertex has degree~4 or 14. Haxell, Seamone, and Verstra\"{e}te~\cite{HSV07} showed that a hamiltonian graph of large enough order $n$ and minimum degree at least $c\log_2 n$, where $c \approx 1.752$, must contain at least two hamiltonian cycles. By combining techniques used in previous sections, we shall now see that there is no threshold at which the minimum degree forces the presence of a superconstant number of hamiltonian cycles. 





\bigskip

\noindent \textbf{Theorem 4.} \emph{For every integer $d \ge 3$ there is an infinite family of hamiltonian $3$-connected graphs with minimum degree~$d$ and with a bounded number of hamiltonian cycles.}

\bigskip

\noindent \emph{Proof.} Let ${\cal G}$ be the infinite family depicted in Fig.~7 and $H \in {\cal G}$. Denote by $H_T$ the graph obtained by applying the replacement operation shown in Fig.~8 to every triangle $T_i$ in $H$.

\begin{center}
\includegraphics[height=55mm]{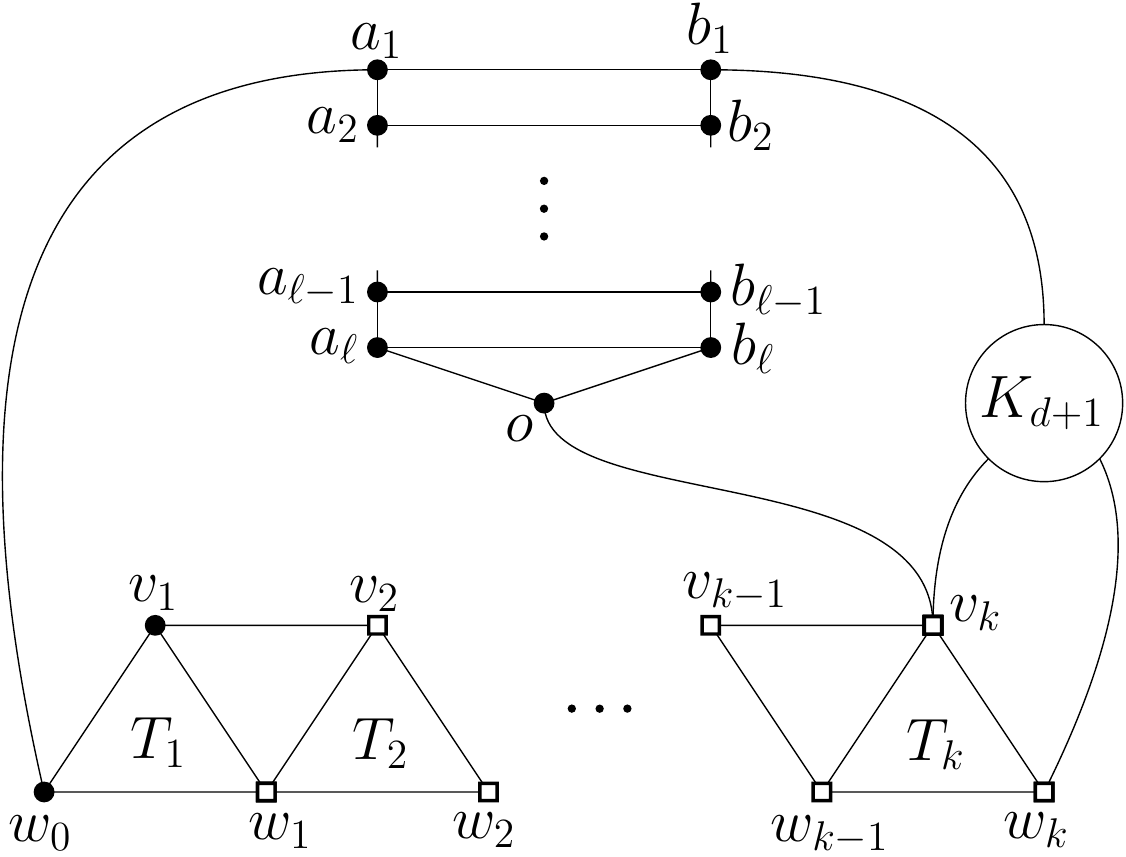}\\
Fig.~7: An infinite family of 3-connected graphs used in the proof of Theorem~4.
\end{center}

\begin{center}
\includegraphics[height=35mm]{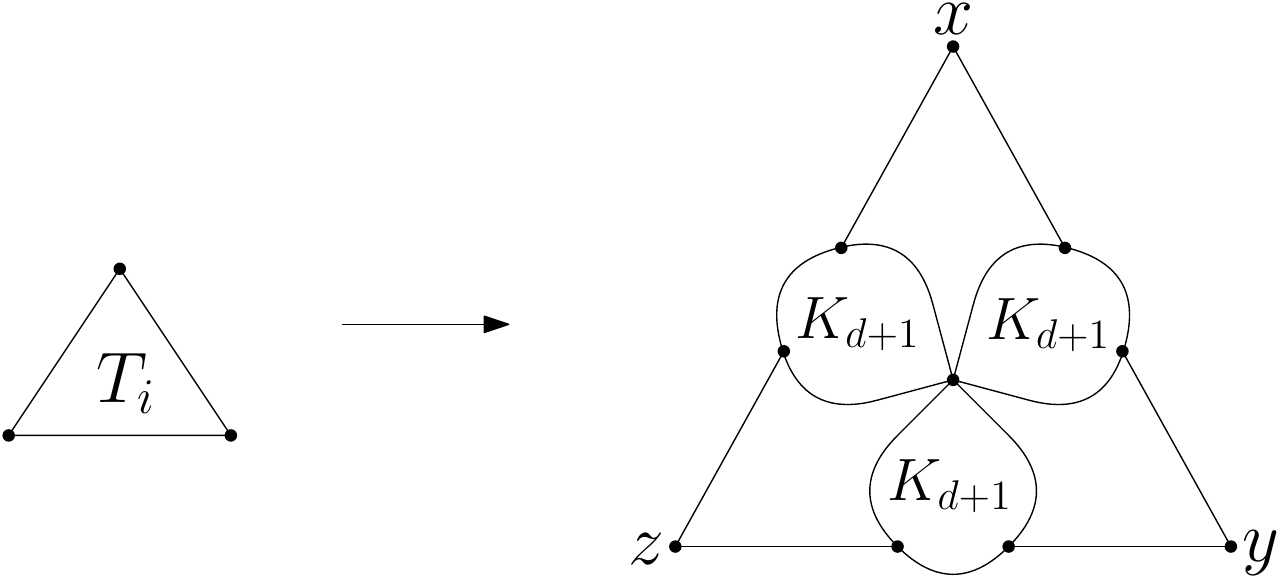}\\
Fig.~8: Replacing the triangle $T_i$ by the graph $T^\Delta$.
\end{center}

We refer to the graph shown on the right-hand side of Fig.~8 as $T^\Delta$ and put $X := \{ x, y, z \}$ where $x,y,z$ are vertices as defined in Fig.~8. We will abbreviate by $(\dagger)$ the observation that in a hamiltonian cycle ${\frak h}$ of a supergraph $S$ of $T^\Delta$ such that among vertices in $T^\Delta$, only $x, y, z$ may have a neighbour in $S - T^\Delta$, the subgraph ${\frak h} \cap T^\Delta$ is a $uv$-path containing $w$, where $u,v,w \in X$ are pairwise distinct. Denote by $W$ ($B$) the set of vertices depicted as white squares (black circular disks) in Fig.~7. For a sufficiently large integer $k$ and any integer $\ell \ge 1$ ($k$ and $\ell$ as introduced in Fig.~7), we can add to $H_T$ edges of the form $vw$ for any $v \in W \cup B$ and $w \in W$, but excluding edges that are already present in $H_T$, such that the resulting graph $G$ has minimum degree $d$ (as $k$ is large enough). In $G$, we denote by $T^\Delta_i$ the copy of $T^\Delta$ which replaced $T_i$. Let ${\frak h}$ be a hamiltonian cycle in $G$. We distinguish between the following three cases---by ($\dagger$), no other cases can occur.

\smallskip

\noindent \textsc{Case 1.} ${\frak h} \cap T^\Delta_1$ is a hamiltonian $v_1w_1$-path. By ($\dagger$), for every $i \in \{ 2, \ldots, k-1 \}$, ${\frak h} \cap T^\Delta_i$ is a hamiltonian $w_{i-1}w_i$-path. Then ${\frak h} \cap T^\Delta_k$ is either a hamiltonian $w_{k-1}v_k$-path or a hamiltonian $w_{k-1}w_k$-path. But then, assuming that the edge $v_1v_k$ or $v_1w_k$ occurs in $G$, $v_1v_k \in E({\frak h})$ or $v_1w_k \in E({\frak h})$, respectively; a contradiction, as $a_1 \notin V({\frak h})$. If neither $v_1v_k$ nor $v_1w_k$ lie in $G$, ${\frak h}$ is not a cycle, again a contradiction.

\smallskip

\noindent \textsc{Case 2.} ${\frak h} \cap T^\Delta_1$ is a hamiltonian $v_1w_0$-path. By ($\dagger$) it is impossible for ${\frak h}$ to span $T^\Delta_2$, a contradiction.

\smallskip

\noindent \textsc{Case 3.} ${\frak h} \cap T^\Delta_1$ is a hamiltonian $w_0w_1$-path. From ($\dagger$) it follows that for every $i \in \{ 1, \ldots, k-1 \}$ we have that ${\frak h} \cap T^\Delta_i$ is a $w_{i-1}w_i$-path spanning $T^\Delta_i$. Using the same arguments as in Case~1 we can infer that $w_0a_1 \in E({\frak h})$. Irrespective of whether ${\frak h} \cap T^\Delta_k$ is a hamiltonian $w_{k-1}v_k$-path or a hamiltonian $w_{k-1}w_k$-path, since $\{ b_1, v_k, w_k \}$ is a 3-cut in $G$, the cycle ${\frak h}$ cannot use the edge $v_kb_\ell$ and must traverse the copy of $K_{d+1}$ located on the right-hand side of Fig.~7, exiting it towards $b_1$. So $a_1b_1 \notin E({\frak h})$. For any integer $\ell \ge 1$, the graph $G[\{ a_1, \ldots, a_\ell, b_1, \ldots, b_\ell, o \}]$ admits exactly one hamiltonian $a_1b_1$-path. It is clear that for a fixed $k$ there exists a constant $c > 0$ such that there are exactly $c$ hamiltonian cycles in $G$. As we can vary $\ell$ whilst adding new edges without altering the number of hamiltonian cycles, the proof is complete. \hfill $\Box$



\section{On 3-regular cyclically 4-edge-connected graphs}

Schwenk~\cite{Sc89} gave a full description of the number of hamiltonian cycles occurring in a generalised Petersen graph $G \in P(n,2)$ ($P(n,k)$ as defined in~\cite{Sc89}); in particular, we have $h(G) = 3$ if and only if $n = 3$ mod.~6. These graphs are cyclically 5-edge-connected, but neither bipartite nor planar. In~\cite{GMZ20} we showed that there exists a planar 3-regular cyclically 4-edge-connected (but non-bipartite) graph on $n$ vertices with exactly four hamiltonian cycles if and only if $n \in \{ 38, 42 \}$ or $n \ge 46$ is even, providing a negative answer to~\cite[Question 1]{CT12} of Chia and Thomassen. A similar result has independently been obtained by Pivotto and Royle~\cite{PR19}. Thomassen~\cite{Th96} described a family of bipartite 3-regular graphs, each containing precisely 16 hamiltonian cycles. These graphs are not planar, and contain a 2-edge-cut; we now give a cyclically 4-edge-connected version of this result.

\bigskip

\noindent \textbf{Proposition 3.} \emph{There exists an infinite family of hamiltonian bipartite $3$-regular cyclically $4$-edge-connected graphs, each containing exactly $16$ hamiltonian cycles.}

\bigskip

\noindent \emph{Proof.} The infinite family ${\cal G}$ of bipartite 3-regular cyclically $4$-edge-connected graphs shown in Fig.~9 is based on a key subgraph of the Ellingham-Horton graph $\Gamma$ on 54 vertices (see~\cite{EH83}); variations thereof have been used in various arguments revolving around small bipartite 3-regular 3-connected graphs that are non-hamiltonian. In Fig.~9, the left-most and right-most part of the graph are to be connected in the obvious way, and an even number of quadrilaterals must occur so that the resulting graph is bipartite.

\begin{center}
\includegraphics[height=45mm]{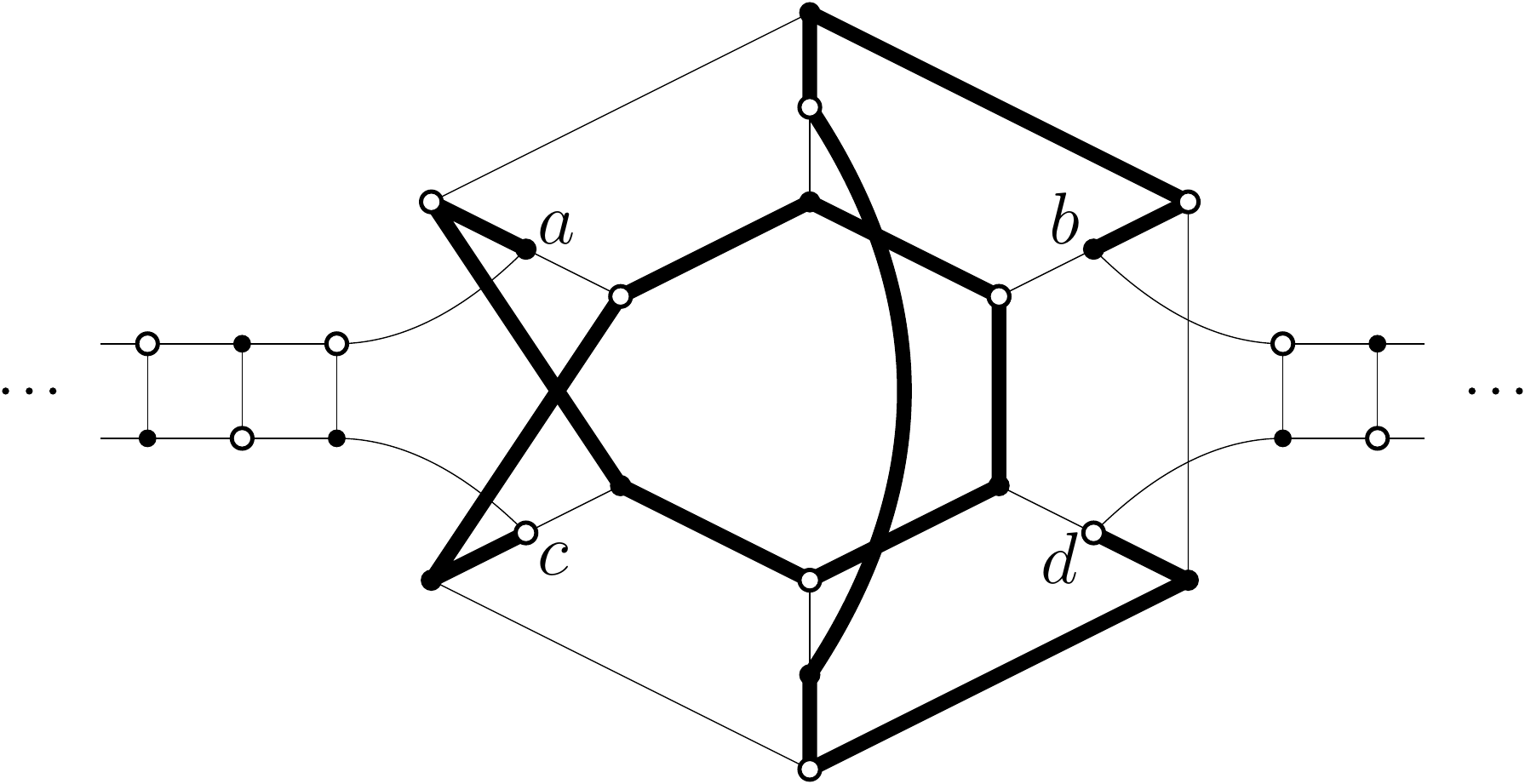}\\
Fig.\ 9: The family ${\cal G}$.
\end{center}

Consider $G \in {\cal G}$. The graph $G - \{ a, b, c, d \}$ has two components (vertices $a,b,c,d$ are as defined in Fig.~9), one of which contains no quadrilateral; we denote it by $C$. Let $H := G[V(C) \cup \{ a, b, c, d \}]$. It follows from the fact that $\Gamma$ is non-hamiltonian that $H$ contains no $ab$-path and $cd$-path whose vertex sets partition $V(H)$, and no $ad$-path and $bc$-path whose vertex sets partition $V(H)$. As shown in Fig.~9, $H$ contains an $ac$-path and a $bd$-path whose vertex sets partition $V(H)$; it is left to the reader to verify that there are exactly 16 such paths. \hfill $\Box$



\section{Regular graphs with a unique longest cycle}

Motivated by work of Chia and Thomassen~\cite{CT12}, we now discuss a different relaxation of Sheehan's conjecture. In contrast to above, where we admitted more than one hamiltonian cycle, we now impose the uniqueness but do not require the cycle to be hamiltonian but simply ask for it to be a longest cycle. By Smith's theorem any hamiltonian 3-regular graph has at least three hamiltonian cycles. But is there a 3-regular graph containing exactly one longest cycle?

Indeed, there is: in~\cite{CT12} it is proven, with a construction based on Petersen's graph, that for each integer $k \ge 2$ there are infinitely many $n$ such that there is a 3-regular graph on $n$ vertices and precisely one longest cycle whose length is $n - k$. We begin by complementing this result with a different solution, shown in Fig.~10, to the above question.

It is straightforward to adapt the constructions from Fig.~10 to $k$-regular graphs with $k \ge 5$ (with a slightly different construction for the even and odd cases) and obtain the following.

\bigskip

\noindent \textbf{Proposition 4.} \emph{For every $k \ge 2$ there exist infinitely many $k$-regular graphs, each containing a unique longest cycle.}

\newpage

\begin{center}
\includegraphics[height=50mm]{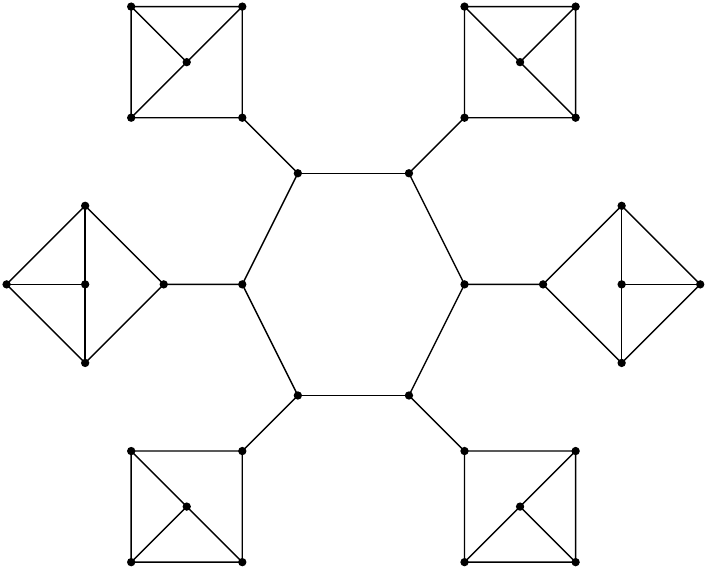} \qquad \includegraphics[height=50mm]{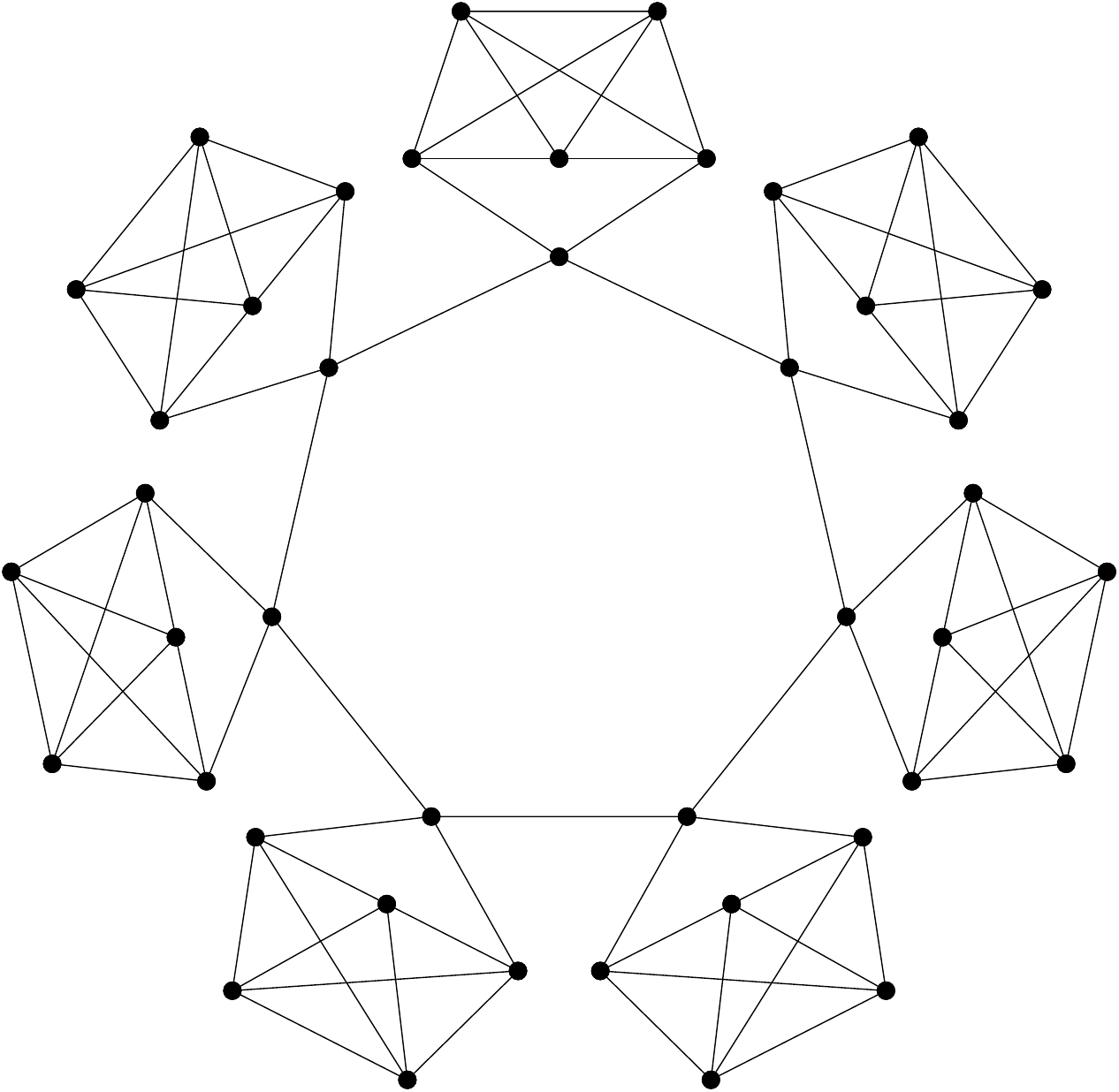}\\[5mm]
\includegraphics[height=16mm]{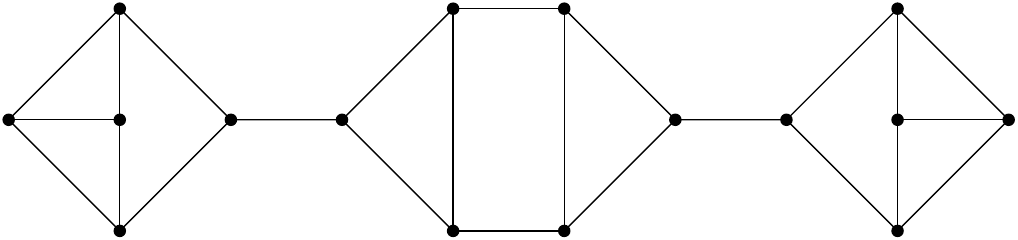}\\
Fig.~10: Top row: 3- and 4-regular graphs with a unique longest cycle, respectively.\\
Bottom: a small 3-regular graph with exactly one longest cycle.
\end{center}


In order to establish a structural property of 3-regular graphs with a unique longest cycle, we require the following lemma. For a graph $G$ we denote by $V_1(G)$ the set of all $1$-valent vertices of $G$, and for a subgraph $H$ in $G$ (allowing $H = G$), we call every vertex in $H$ which has degree~1 in $H$ an \emph{$H$-leaf}.

\bigskip

\noindent \textbf{Lemma 3.} \emph{A graph $F$ with maximum degree at most~$3$ and $|V_1(F)| \ge 2$ contains $\lfloor |V_1(F)|/2 \rfloor$ pairwise disjoint paths all of whose end-vertices lie in $V_1(F)$.}

\bigskip

\noindent \emph{Proof.} For graphs with maximum degree at most~$3$ and exactly two 1-valent vertices the statement is clearly true. Assume the statement to be true for graphs with maximum degree at most~$3$ and at most $\ell - 1$ 1-valent vertices, and let $G$ be a graph with maximum degree at most 3 and $|V_1(G)| = \ell \ge 3$. Let $T$ be a spanning tree of $G$. Clearly, every $G$-leaf is a $T$-leaf, but not every $T$-leaf needs to be a $G$-leaf. In $T$, consider the unique path $v_1 \ldots v_k$ between distinct $G$-leaves $v_1$ and $v_k$. For $v_i$ with $\deg_T(v_i) = 3$, denote by $\ell(i)$ the number of $G$-leaves in the unique $\{ v_i \}$-fragment $F_i$ of $T$ containing neither $v_1$ nor $v_k$. Define the function $f : \{ 2, \ldots, k-1 \} \rightarrow \mathbb{N}$ with
$$f(i) := \begin{cases}
            0 & {\rm if} \ \deg_T(v_i) = 2,\\
            \ell(i) & {\rm else.}
          \end{cases}$$
Our aim is to find $G$-leaves $v,w$ such that for the $vw$-path $P$ in $T$ we have that either all components of $T - P$ or all but one component of $T - P$ contain(s) an even number of $G$-leaves. This implies that either all components of $G - P$ or all but one component of $G - P$ contain(s) an even number of $G$-leaves, from which the lemma's statement follows by applying the assumption made in this proof's second sentence to each component.

If no $f(i)$ is odd then we are done, as we can choose $v = v_1$ and $w = v_k$. So assume at least one $f(i)$ to be odd. Consider the path $Q'$ from $v_1$ to $v_i$ with $i$ minimal such that $f(i)$ is odd. $F_i$ (which is a tree) contains at most $\ell - 2$ $G$-leaves as $T$ contains exactly $\ell$ $G$-leaves and neither $v_1$ nor $v_k$ lie in $F_i$. The vertex $v_i$ is not a $G$-leaf but it is an $F_i$-leaf. Let $S$ be the union of $\{ v_i \}$ and the set of all $G$-leaves in $F_i$. Since $f(i)$ is odd, $|S|$ is even. We know that there exist $|S|/2$ pairwise disjoint paths in $F_i$ all of whose end-vertices lie in $S$. Every $G$-leaf in $F_i$ is contained in one of these paths. Exactly one of these paths, say $Q''$, has $v_i$ as endpoint. Let the other endpoint of $Q''$ be $v'$. We have $Q' \cap Q'' = \{ v_i \}$.

Put $Q := Q' \cup Q''$. Setting $v = v_1$ and $w = v'$ we are done, as either all components or all but one component of $T - Q$ contain(s) an even number of $G$-leaves: $f(j)$ is even for all $j < i$, and by the choice of $w$, for every $u \in V(Q'') \setminus \{ v_i, w \}$ with $\deg_T(u) = 3$ the unique $\{ u \}$-fragment of $T$ containing neither $v$ nor $w$ must contain an even number of $G$-leaves. \hfill $\Box$

\bigskip

\noindent \textbf{Theorem 5.} \emph{Let $G$ be a $3$-regular graph. If $G$ has a unique longest cycle ${\frak c}$, then at least two components of $G - E({\frak c})$ have an odd number of vertices in ${\frak c}$. There exist infinitely many $3$-regular graphs with a unique longest cycle and exactly two such components.}

\bigskip

\noindent \emph{Proof.} Assume there exists a 3-regular graph $G$ containing a unique longest cycle ${\frak c}$ in which every component of $G - E({\frak c})$ has an even number of vertices in ${\frak c}$ (reductio ad absurdum). The possibly disconnected graph $G - E({\frak c})$ consists of a set of graphs in which every vertex has degree either 1 or 3, and in which there are at least two vertices of degree~1 by hypothesis. In this set, we ignore all graphs isomorphic to $K_2$ and consider the remaining elements, which we call $F_1, \ldots, F_k$. Note that $k \ge 1$ as ${\frak c}$ is non-hamiltonian by Smith's theorem stating that 3-regular graphs cannot be uniquely hamiltonian.

We consider for each $F_i$ a set ${\cal P}_i$ of $\lfloor |V_1(F_i)|/2 \rfloor$ pairwise disjoint paths with the properties described in Lemma~3. Every path in ${\cal P}_i$ has at least three vertices as $F_i \ne K_2$. We obtain the graph $G'$ by removing from $G$ all vertices which are not in ${\frak c}$ or in a path from $\bigcup_i {\cal P}_i$, and by replacing, in each $F_i$, every path $P \in {\cal P}_i$ by an edge between the endpoints of $P$ (which lie in ${\frak c}$). The graph $G'$ is 3-regular.

${\frak c}$ is still a unique longest cycle in $G'$: assume there would be a cycle ${\frak c}' \ne {\frak c}$ in $G'$ of length at least $|V({\frak c})|$. In ${\frak c}'$, consider every edge that replaced a path in the above argument (there must be at least one such edge since otherwise ${\frak c}'$ would have occurred already in $G$) and revert the replacement, i.e.\ replace the edge by the path it replaced. We thus obtain a cycle ${\frak c}'' \ne {\frak c}$ in $G$ of length at least $|V({\frak c})|$, a contradiction. Moreover, due to the replacements we performed, ${\frak c}$ is a hamiltonian cycle in $G'$, so $G'$ is a 3-regular graph containing exactly one hamiltonian cycle; but this is impossible by Smith's theorem.

Now suppose there exists a 3-regular graph $G$ containing a unique longest cycle ${\frak c}$ in which exactly one component $F$ of $G - E({\frak c})$ has an odd number of vertices in ${\frak c}$ (reductio ad absurdum). We proceed exactly as above for all components with an even number of vertices in ${\frak c}$, transforming them into disjoint unions of $K_2$'s. By Lemma~3, there exists a vertex $x \in V_1(F) \cap V({\frak c})$ such that $F$ contains $(|V_1(F)| - 1)/2$ pairwise disjoint paths whose endpoints span $V_1(F) \setminus \{ x \}$. We note that none of these paths contain $x$, and that $(|V_1(F)| - 1)/2 = 0$ is possible. We delete all vertices that are not in ${\frak c}$ or contained in one of the aforementioned paths, and contract each path to $K_2$. After performing these operations, we add a new vertex $y$, the edge $xy$ and a loop from $y$ to itself. We obtain a 3-regular multigraph $G''$ with a unique longest cycle of length $|V(G'')| - 1$. But Thomason showed in~\cite{Th78} that if a 3-regular multigraph contains a vertex-deleted subgraph which has an odd number of hamiltonian cycles, then the graph itself must be hamiltonian. Thus, we have obtained a contradiction.

That there exist 3-regular graphs containing a unique longest cycle ${\frak c}$ such that exactly two components of $G - E({\frak c})$ have an odd number of vertices in ${\frak c}$ is shown in Fig.~10, and it is straightforward to deduce from that construction the existence of an infinite family with the stated properties. Chia and Thomassen~\cite{CT12} give in their paper a structurally very different such example: consider Petersen's graph $P$ and a vertex $x$ in $P$. Inflate every vertex in $P$ except $x$ to a triangle so as to obtain $H$, a graph containing exactly two longest cycles. In $H$, consider an edge $vw$ which lies on one longest cycle of $H$ but not the other. Let $H_1$ and $H_2$ be disjoint copies of $H - vw$ such that $v_i$ ($w_i$; $x_i$) is the copy of $v$ ($w$; $x$) in $H_i$. Then the graph $G := H_1 \cup H_2 + v_1v_2 + w_1w_2$ is 3-regular, contains a unique longest cycle ${\frak c}$, and the components of $G - E({\frak c})$ are all isomorphic to $K_2$ with two exceptions, which are isomorphic to $K_{1,3}$; the latter two graphs contain the vertices $x_1$ and $x_2$. \hfill $\Box$

\bigskip

The 3-regular graphs with a unique longest cycle presented above have connectivity~1, while every member of the infinite family given by Chia and Thomassen~\cite{CT12} has connectivity~2. We now show an analogous result for connectivity~3, despite the dramatic (but, as explained in~\cite{AT97}, perhaps to be expected) difference between the minimum number of cycles present in 3-regular 2-connected graphs and the minimum number of cycles in 3-regular 3-connected graphs---the former is quadratic while the latter is superpolynomial in the graphs' order.

\bigskip

\noindent \textbf{Proposition 5.} \emph{There exist infinitely many integers $n$ such that there exists a $3$-regular $3$-connected graph of order~$n$ and containing a unique longest cycle whose length is $n - 2$.}

\bigskip

\noindent \emph{Proof.} Our proof begins in the same way the proof of Theorem~5 ended: consider Petersen's graph $P$ and distinct non-adjacent vertices $x,y \in V(P)$. We obtain the graph $H$ by replacing every vertex in $P$ except $x$ and $y$ by a triangle and removing $x$. One deduces from the properties of the Petersen graph that one can label the neighbours of $x$ by $x_1, x_2, x_3$ such that (i) in $H$, there is no hamiltonian $x_ix_j$-path for any $i,j \in \{ 1, 2, 3 \}$; (ii) in $H - y$ there is exactly one hamiltonian $x_1x_2$-path, exactly one hamiltonian $x_1x_3$-path, and no hamiltonian $x_2x_3$-path; and (iii) there exists for no $v \in V(H) \setminus \{ y \}$ a hamiltonian $x_ix_j$-path in $H - v$, for any $i,j \in \{ 1, 2, 3 \}$.

Consider two disjoint copies of $H$ which we call $H'$ and $H''$. For a vertex $v$ in $H$, denote by $v'$ ($v''$) the corresponding vertex in $H'$ ($H''$). Put $M := \{ x'_1x''_2, x'_2x''_1, x'_3x''_3 \}$. The graph $G := H' \cup H'' + M$ is clearly 3-regular and 3-connected. Let ${\frak c}$ be a longest cycle in $G$. By (i) and (ii), the length of ${\frak c}$ is 48 (while the order of $G$ is 50), so ${\frak c}$ must traverse $M$. It does so by using exactly two of $M$'s edges. If $x'_3x''_3 \in E({\frak c})$, then by (ii), ${\frak c} \cap H'$ is an $x'_1x'_3$-path. But by (ii) and (iii), there is no $x''_2x''_3$-path in $H''$ of length at least 24, so we can conclude that $x'_3x''_3 \notin E({\frak c})$, whence, $x'_1x''_2, x'_2x''_1 \in E({\frak c})$. The uniqueness of ${\frak c}$ now follows from (ii).

From this example one obtains an infinite family by successively considering vertices lying on the longest cycle and inflating them to triangles. \hfill $\Box$

\bigskip

By the already mentioned results of Smith and Thomason, the ``$n - 2$'' in Proposition~5 is best possible.

\section{Discussion}

\noindent \textbf{1.} By \emph{triangulation} we mean a plane graph in which every face is a triangle. Hakimi, Schmeichel, and Thomassen~\cite{HST79} proved that there exist infinitely many triangulations, each containing exactly four hamiltonian cycles (while Kratochvil and Zeps~\cite{KZ88} proved that besides $K_3$ and $K_4$, hamiltonian triangulations cannot have fewer than four hamiltonian cycles). Each of these triangulations $T$ contains exactly six cubic vertices. We consider $T$ to be embedded in the plane. By replacing in $T$ each triangle together with the cubic vertex it contains by an octahedron, it follows that there is a positive constant $c$ such that there are infinitely many triangulations of minimum degree~$4$, each containing exactly $c$ hamiltonian cycles. There is however no obvious way in which their result can be adapted in order to show that the same conclusion holds for minimum degree~5 (since 4-valent vertices play a pivotal structural role). We observe that, by using the techniques introduced in this manuscript, one can prove that there exists an infinite family of hamiltonian triangulations of minimum degree~$5$ with a bounded number of hamiltonian cycles.

\smallskip



\noindent \textbf{2.} The results in Section~2 suggest that between 4- and 5-regularity, a dramatic shift occurs with respect to the way the hamiltonian graph's order relates to its minimum number of hamiltonian cycles. Is there an infinite family of hamiltonian graphs of minimum degree at least~4 and maximum degree~5, and with a bounded number of hamiltonian cycles, in which, asymptotically, the number of non-4-valent vertices does not vanish?

\smallskip

\noindent \textbf{3.} A famous result of Tutte states that planar 4-connected graphs are hamiltonian. Thomassen extended this by proving that every planar graph with minimum degree at least~4 in which every vertex-deleted subgraph is hamiltonian, must itself be hamiltonian~\cite{Th78-2}. In~\cite{Za}, we showed that such graphs must contain at least three hamiltonian cycles. For a long time, only a constant number of hamiltonian cycles was guaranteed to exist in planar 4-connected graphs; however, Brinkmann and Van Cleemput recently proved that planar 4-connected graphs contain at least a linear number of hamiltonian cycles~\cite{BV21}. Does this result extend to planar graphs with minimum degree at least~4 in which every vertex-deleted subgraph is hamiltonian---or are there such graphs with ``few'' hamiltonian cycles?

\smallskip

\noindent \textbf{4.} In problems such as the ones discussed in this article it is common to investigate girth restrictions, one example being Cantoni's conjecture that every planar 3-regular graph with exactly three hamiltonian cycles contains a triangle (see~\cite{Tu76}), which in~\cite{GMZ20} was confirmed to hold for graphs of order at most~48. In view of Haythorpe's conjecture on 4-regular graphs, it would be of interest to determine the minimum number of hamiltonian cycles in hamiltonian 4-regular triangle-free graphs. (In~\cite{GMZ20} these numbers were computed for all graphs up to order~21.) We know by a theorem of Thomassen~\cite{Th96} that every hamiltonian \emph{bipartite} graph of minimum degree at least 4 and girth $g$ has at least $(3/2)^{g/8}$ hamiltonian cycles. Finally, in the light of Section~6, we ask: are there 3-regular 2-connected triangle-free graphs containing a unique longest cycle?

\vspace{5mm}

\noindent \textbf{Acknowledgements.} Thanks are due to the referees for their helpful comments. The research presented in this paper was supported by a Postdoctoral Fellowship of the Research Foundation Flanders (FWO).

\newpage

\section{Appendix}

\begin{center}
\includegraphics[height=75mm]{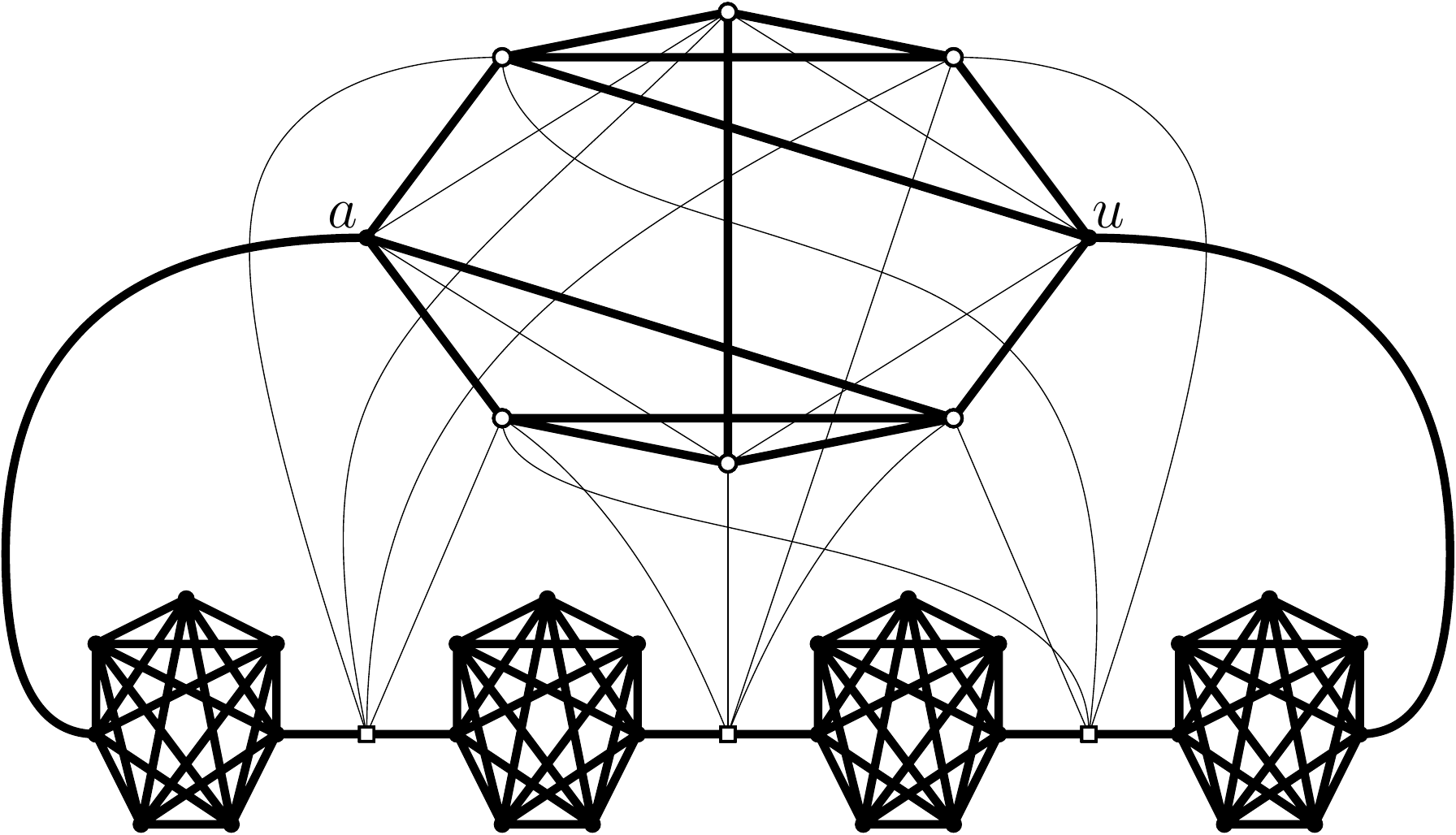}\\
Fig.~A: A 6-regular graph with exactly $5 \cdot (5!)^4$ hamiltonian cycles. This graph is used in the proof of Theorem~3.
\end{center}

\vspace{1cm}

\begin{center}
\includegraphics[height=70mm]{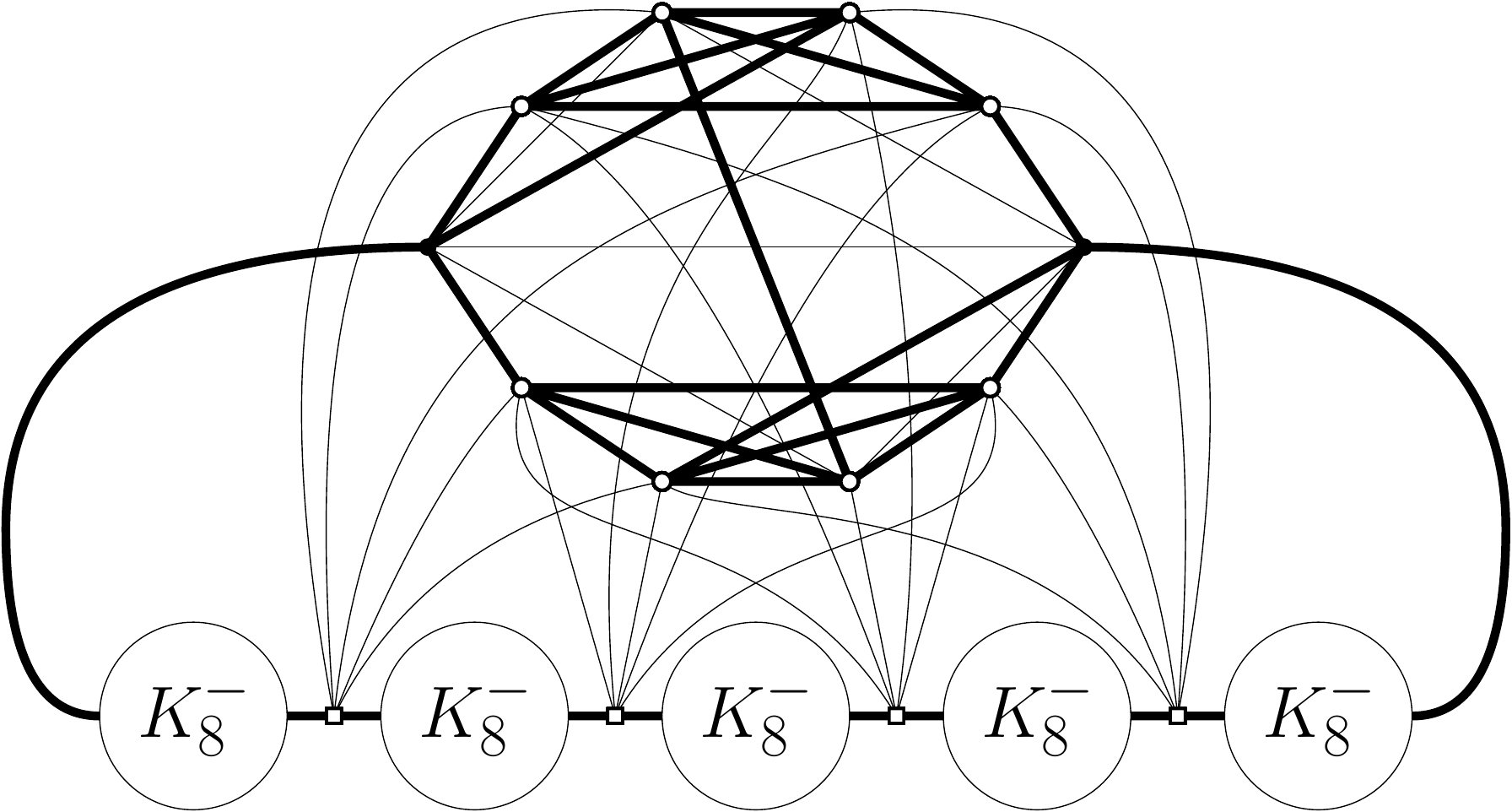}\\
Fig.~B: A 7-regular graph $G$ with exactly $16 \cdot (6!)^5$ hamiltonian cycles; $K_8^-$ stands for a complete graph on eight vertices minus the unique (horizontal) edge whose removal renders $G$ 7-regular. This graph is used in the proof of Theorem~3.
\end{center}


\begin{thebibliography}{99}

\bibitem{AT97}
R. E. L. Aldred and C. Thomassen.
On the Number of Cycles in 3-Connected Cubic Graphs.
\emph{J. Combin. Theory, Ser. B} \textbf{71} (1997) 79--84.

\bibitem{BS20}
R. D. Barish and A. Suyama.
Counting Hamiltonian Cycles on Quartic 4-Vertex-Connected Planar Graphs.
\emph{Graphs Combin.} \textbf{36} (2020) 387--400.

\bibitem{Bo95}
J. A. Bondy.
Basic graph theory.
In: Handbook of Combinatorics (Eds.:~M.~Gr\"otschel, L. Lov\'asz, R. L. Graham), pp.~3--110, Amsterdam, North-Holland, 1995.

\bibitem{BJ98}
J. A. Bondy and B. Jackson.
Vertices of Small Degree in Uniquely Hamiltonian Graphs.
\emph{J. Combin. Theory, Ser. B} \textbf{74} (1998) 265--275.


\bibitem{BV21}
G. Brinkmann and N. Van Cleemput.
4-connected polyhedra have at least a linear number of hamiltonian cycles.
\emph{Europ. J. Combin.} \textbf{97} (2021) Aricle number 103395.


\bibitem{CT12}
G. L. Chia and C. Thomassen.
On the number of longest and almost longest cycles in cubic graphs.
\emph{Ars Combin.} \textbf{104} (2012) 307--320.

\bibitem{EH83}
M. N. Ellingham and J. D. Horton.
Non-hamiltonian 3-connected cubic bipartite graphs.
\emph{J. Combin. Theory, Ser. B} \textbf{34} (1983) 350--353.

\bibitem{ES80}
R. C. Entringer and H. Swart.
Spanning cycles of nearly cubic graphs.
\emph{J. Combin. Theory, Ser. B} \textbf{29} (1980) 303--309.

\bibitem{Fl14}
H. Fleischner.
Uniquely Hamiltonian Graphs of Minimum Degree 4.
\emph{J. Graph Theory} \textbf{75} (2014) 167--177.

\bibitem{GKN19}
A. Gir\~{a}o, T. Kittipassorn, and B. Narayanan.
Long cycles in Hamiltonian graphs.
\emph{Israel J. Math.} \textbf{229} (2019) 269--285.

\bibitem{GMZ20}
J. Goedgebeur, B. Meersman, and C. T. Zamfirescu.
Graphs with few hamiltonian cycles.
\emph{Math. Comp.} \textbf{89} (2020) 965--991.

\bibitem{HST79}
S. L. Hakimi, E. F. Schmeichel, and C. Thomassen.
On the Number of Hamiltonian Cycles in a Maximal Planar Graph.
\emph{J. Graph Theory} \textbf{3} (1979) 365--370.

\bibitem{HSV07}
P. Haxell, B. Seamone, and J. Verstraete.
Independent dominating sets and hamiltonian cycles.
\emph{J. Graph Theory} \textbf{54} (2007) 233--244.

\bibitem{Ha18}
M. Haythorpe.
On the Minimum Number of Hamiltonian Cycles in Regular Graphs.
\emph{Experim. Math.} \textbf{27} (2018) 426--430.

\bibitem{KZ88}
J. Kratochvil and D. Zeps.
On the number of Hamiltonian cycles in triangulations.
\emph{J. Graph Theory} \textbf{12} (1988) 191--194.

\bibitem{Me73}
G. H. J. Meredith.
Regular $n$-Valent $n$-Connected NonHamiltonian Non-$n$-Edge-Colorable Graphs.
\emph{J. Combin. Theory, Ser. B} \textbf{14} (1973) 55--60.


\bibitem{PR19}
I. Pivotto and G. Royle.
Highly-connected planar cubic graphs with few or many Hamilton cycles.
\emph{Discrete Math.} \textbf{342} (2019) 111608.

\bibitem{Ro12}
G. Royle. Two problems on Hamilton cycles. Blog ``SymOmega'', 1 Feb.~2012.\\
\url{https://symomega.wordpress.com/2012/02/01/two-problems-on-hamilton-cycles}

\bibitem{Sc89}
A. J. Schwenk.
Enumeration of Hamiltonian cycles in certain generalized Petersen graphs.
\emph{J. Combin. Theory, Ser. B} \textbf{47} (1989) 53--59.

\bibitem{Sh75}
J. Sheehan.
The multiplicity of Hamiltonian circuits in a graph.
In: Recent advances in graph theory (Proc. Second Czechoslovak Sympos., Prague, 1974), Academia, Prague (1975) 477--480.

\bibitem{Th78}
A. G. Thomason.
Hamiltonian Cycles and Uniquely Edge Colourable Graphs.
\emph{Ann. Discrete Math.} \textbf{3} (1978) 259--268.

\bibitem{Th78-2}
C. Thomassen.
Hypohamiltonian graphs and digraphs.
In: Theory and Applications of Graphs, Lecture Notes in Mathematics \textbf{642}, Springer, Berlin (1978) 557--571.

\bibitem{Th96}
C. Thomassen.
On the Number of Hamiltonian Cycles in Bipartite Graphs.
\emph{Combin. Probab. Comput.} \textbf{5} (1996) 437--442.

\bibitem{Th98}
C. Thomassen.
Independent Dominating Sets and a Second Hamiltonian Cycle in Regular Graphs.
\emph{J. Combin. Theory, Ser. B} \textbf{72} (1998) 104--109.

\bibitem{TZ}
C. Thomassen and C. T. Zamfirescu.
4-regular 4-connected Hamiltonian graphs with a bounded number of Hamiltonian cycles.
\emph{Australasian J. Combin.} \textbf{81} (2021) 334--338.

\bibitem{Tu46}
W. T. Tutte.
On Hamiltonian circuits.
\emph{J. London Math. Soc.} \textbf{21} (1946) 98--101.

\bibitem{Tu76}
W. T. Tutte.
Hamiltonian circuits.
In: \emph{Colloquio Internazionale sulle Teorie Combinatorie} (Roma, 1973), Tomo I, Atti dei Convegni Lincei, No.~17, Accad.\ Naz.\ Lincei, Rome (1976) 193--199.

\bibitem{VZ18}
N. Van Cleemput and C. T. Zamfirescu.
Regular non-hamiltonian polyhedral graphs.
\emph{Appl. Math. Comput.} \textbf{338} (2018) 192--206.

\bibitem{Wa13}
A. Wagner.
On the Existence of a Second Hamilton Cycle in Hamiltonian Graphs With Symmetry.
M.Sc.~Thesis, University of Ottawa, Canada, 2013.

\bibitem{Za}
C. T. Zamfirescu.
On the hamiltonicity of a planar graph and its vertex-deleted subgraphs.
Submitted.

\end{thebibliography}
\end{document}